\documentclass[11pt]{article}

\usepackage{amsmath,amssymb,amsfonts}
\usepackage{mathtools}

\usepackage{amsthm}
\newtheorem{theorem}{Theorem}[section]
\newtheorem{lemma}[theorem]{Lemma}

\newtheorem{definition}[theorem]{Definition}
\newtheorem{conjecture}[theorem]{Conjecture}

\newtheorem*{remark}{Remark}

\usepackage{listings}
\usepackage{xcolor}
\usepackage{graphicx}
\usepackage{float}

\usepackage{booktabs}
\usepackage{array}

\usepackage{caption}
\usepackage[utf8]{inputenc}
\usepackage[T1]{fontenc}

\usepackage{geometry}
\usepackage{microtype}

\usepackage{url}
\usepackage{hyperref}
\hypersetup{
    colorlinks=true,
    linkcolor=blue,
    citecolor=blue,
    urlcolor=blue,
    breaklinks=true,
    pdftitle={Combinatorial and Gaussian Foundations of Rational Nth Root Approximations},
    pdfauthor={Isaac Wolford},
    pdfsubject={Mathematics, Numerical Analysis}
}

\usepackage{cleveref}
\crefname{figure}{Figure}{Figures}
\crefname{table}{Table}{Tables}
\crefname{equation}{Equation}{Equations}
\crefname{theorem}{Theorem}{Theorems}
\crefname{lemma}{Lemma}{Lemmas}
\crefname{definition}{Definition}{Definitions}
\crefname{conjecture}{Conjecture}{Conjectures}

\title{Combinatorial and Gaussian Foundations of Rational Nth Root Approximations: Theorems and Conjectures}
\author{Isaac Wolford\\
\texttt{isaacmw@southern.edu}\\
\small ORCID: \href{https://orcid.org/0009-0003-4231-0912}{0009-0003-4231-0912}}
\date{\today}

\begin{document}

\maketitle
\begin{abstract}
We present an approach (the biroot method) for nth root approximation that yields closed-form rational functions with coefficients derived from binomial structures, Gaussian functions, or qualifying DAG structures. The method emerges from an analysis of Newton's method applied to square root extraction, revealing that successive iterations generate coefficients following rows of Pascal's triangle in an alternating numerator-denominator pattern. After further exploration of these patterns, we formulate three main conjectures: (1) the Binomial Biroot Conjecture establishing the fundamental alternating coefficient structure to approximate nth roots (for which we prove the square root case and optimal parameter conditions), (2) the Gaussian Biroot Conjecture, and (3) the General DAG Biroot Conjecture showing a structural invariance to nth root approximation using arbitrary linearly-constructed directed acyclic graphs (DAGs). Computational evidence demonstrates superior convergence properties when compared to Taylor series and Padé approximations, especially considering the more direct and less computationally intensive approach to the biroot function construction. A computational framework has been developed to support systematic exploration of the biroot method's parameter space and to enable extensive numerical and symbolic analysis. The method provides both theoretical insights and computational significance by connecting combinatorial structures to nth root rational approximation theory.
\end{abstract}


\section{Introduction}
Classical methods for approximating nth roots—Taylor series, Padé approximations, and Newton's method—each possess inherent limitations. Taylor series exhibit poor convergence beyond narrow intervals, Padé approximations (including Chebyshev variants) require solving complex systems of equations, and Newton's method provides only recursive solutions without closed-form expressions.

This paper introduces the biroot method, an nth root approximation approach that addresses these limitations through closed-form rational functions with superior convergence properties. The method emerges from a systematic analysis of Newton's method applied to square root extraction, revealing that successive symbolic iterations generate coefficients following rows of Pascal's triangle (PT for short) in an alternating numerator-denominator pattern. Generalizing these patterns produces closed-form expressions for nth root approximations. We prove convergence for the square root case, and propose a proof pathway for the nth root case. We also prove the optimal parameter conditions for the binomial biroot. The Central Limit Theorem then provides theoretical justification for using continuous Gaussian functions to sample coefficient values, yielding improved convergence properties. Furthermore, generalizing the approach to arbitrary linearly constructed PT-like DAG structures exhibits notable invariance properties for nth root approximation.

Extensive computational evidence supports these insights. Systematic exploration across multiple parameter spaces demonstrates superior convergence rates compared to classical methods, with particularly strong performance for the Gaussian variant. This computational investigation also enables us to formulate conjectures regarding optimal conditions, invariance properties, and equivalence to Padé approximations. Our analysis framework\footnote{Available in the supplementary sources: https://doi.org/10.5281/zenodo.16878055} facilitates comprehensive comparison with existing methods and provides tools for error analysis and convergence studies through a cohesive interface.

Beyond computational validation, our work suggests broader theoretical insights regarding the fundamental structure underlying nth root approximation methods. Such methods as Newton's method, Padé approximations, and Taylor series may be understood as discrete manifestations of a deeper continuous framework. In this view, the Gaussian biroot may represent the theoretical limit of optimal coefficient distribution for rational approximation of radical functions.

Our investigation establishes three conjectures: (1) the Binomial Biroot Conjecture formalizing convergence of the basic alternating structure, (2) the Gaussian Biroot Conjecture establishing seemingly optimal continuous coefficient distributions, and (3) the DAG Biroot Conjecture extending to arbitrary linearly-constructed graph frameworks. The remainder of this paper develops these ideas systematically.

\begin{definition}[Biroot Approximant]
A rational function of the form $\frac{P_m(x)}{Q_{k}(x)}, k \in \{m, m-1\}$ where the coefficients of polynomials $P_m$ and $Q_k$ are derived from combinatorial structures (binomial coefficients, Gaussian samples, or DAG node values) arranged in an alternating fashion between numerator and denominator, to approximate $x^{1/n}$. We use the term \textit{biroot} to refer to such approximants, emphasizing the binomial origin, alternating coefficient pattern, and root approximation purpose. For the sake of simplicity, $\beta_m^n(x, c)$ is used to indicate a biroot approximation of index $n$, order $m$, and expansion point $c$. In general, $c$ is $1$ by default, $n$ is $2$ by default, and $m$ can be inferred from the context. The formula for $\beta$ can also be inferred from the section: Binomial, Gaussian, or DAG.
\end{definition}

\section{Background and Related Work}
Classical approaches to nth root approximation each possess distinct advantages and limitations that inform the development of the biroot method. This section examines both established methods—Newton's method, Taylor polynomials, and Padé approximations—and recent developments in rational approximation theory.

\subsection{Newton's method}
Newton's method provides an efficient approach for nth root computation \cite{Galantai2000,Dubeau1998,Dubeau2009}. The general iteration formula is as follows.
\[x_{k+1}=x_k - \frac{f(x_k)}{f'(x_k)}\]

Applying this to find $\sqrt[n]{a}$ using $f(x) = x^n - a$ gives the general form of Newton's method applied to nth roots:
\begin{align*}
x_{k+1} = \frac{(n-1)x_k^n + a}{nx^{n-1}_k}
\end{align*}

While Newton's method converges rapidly, it provides only recursive solutions. This recursive nature obfuscates symbolic patterns and requires iterative computation for each evaluation, limiting its utility in applications requiring closed-form expressions.

\subsection{Taylor Polynomials}
\cite{Taylor2021,Magyar2025}
Consider the nth root function $x^{1/n}$. A general expression of the Taylor polynomial to the $m$th degree for this function is as follows.
\begin{align*}
T(x) \sim  x^{1/n} = \sum^m_{k=0} \frac{\left(\frac{1}{n}\right)^{\underline{k}} \cdot c^{1/n-k}}{k!} (x-c)^k
\end{align*}
where $\left(\frac{1}{n}\right)^{\underline{k}} \cdot x^{\frac{1}{n}-k}$ is the $k$th derivative of $x^{1/n}$. 

It should be noted that it is possible for a Taylor series to exist that approximates the nth root at any fixed point $x>0$. This means that it is possible to construct Taylor polynomials that converge over large intervals. Computing these polynomials, however, becomes increasingly unwieldy as greater fixed points $x$ need to be approximated.

For instance, to approximate $\sqrt{10}$ with accuracy $10^{-6}$ using Taylor series centered at $c=1$, one would need impractically high degree polynomials since $10 \notin (0,2c)$. Even expanding at $c=5$ to include the target value requires computing numerous high-order derivatives and results in unwieldy polynomial expressions. This limitation motivates our search for closed-form rational approximations that maintain accuracy across extended intervals without requiring derivative computations or high polynomial degrees.

\subsection{Padé Approximations}
Padé approximations are rational function approximations with superior convergence compared to Taylor series \cite{Chen1989,Abrate2013,Baker1996,Fornberg}. A Padé approximant $P_{[p/q]}^{(c)}$ has numerator degree $p$, denominator degree $q$, and expansion point $c$. The construction process requires solving a system of equations. Given polynomials $P$ and $Q$ with unknown coefficients, a Taylor polynomial $A$ of order $p+q$ is constructed, and coefficients are determined by solving for the coefficients in the equation $P = QA$. The solved coefficients are then used to construct the rational function \(P/Q\). While this yields excellent rational approximations, the steps required to construct such approximants are not straight forward and cannot be represented by a single, elegant, closed-form formula.

Notably, Padé approximations appear to exhibit a direct connection to our biroot method. For example, the $[2/2]$ Padé approximant for $\sqrt{x}$ expanded at $1$ yields:
\[\sqrt{x} \approx \frac{1 + 10x + 5x^2}{5 + 10x + x^2}\]
This expression is identical to a specific biroot approximation, suggesting deeper theoretical connections explored in our Padé Equivalence Conjecture. This relationship indicates that biroot methods may provide direct construction of certain Padé approximants without requiring system solutions.

\subsection{Recent Developments in Closed-Form Root Approximation}
Recent work by Shunia \cite{shunia2024} has introduced an approach to nth root computation using Kronecker substitution and polynomial quotient rings, representing noteworthy progress toward closed-form root expressions. Shunia's method derives a limit formula for $\sqrt[n]{a}$ characterized by modular exponentiations
\[\sqrt[n]{a} = \lim_{k \to \infty} \frac{(k^{kn} + 1)^{kn+1} \bmod (k^{kn^2} - a)}{(k^{kn} + 1)^{kn} \bmod (k^{kn^2} - a)} - 1\]
and, more relevantly to our work, conjectures a direct closed-form expression for the integer part of nth roots
\[\lfloor\sqrt[n]{a}\rfloor = \left\lfloor \frac{(a^{2an} + 1)^{2an+1} \bmod (a^{2an^2} - a)}{(a^{2an} + 1)^{2an} \bmod (a^{2an^2} - a)} - 1 \right\rfloor\]
This conjecture represents a breakthrough in closed-form root computation, achieving exact integer results through modular arithmetic. Shunia's approach differs from our method in several key aspects: his technique employs modular exponentiations with bases growing as $a^{2an}$, focuses specifically on integer parts of roots, and provides a single formula per root index.

In contrast, our biroot method offers a complementary approach through rational approximations constructed from combinatorial patterns discovered in Newton's method iterations and Padé approximants, and extends this pattern to continuous Gaussian functions and generalized PT-like DAG structures. We also establish optimal expansion conditions through combinatorics. Our approach offers a broader framework for understanding and utilizing the combinatorial foundations underlying nth root approximation.

\section{Conjectures and Theorems}
The first conjecture (Binomial Biroot) consists of the initial reasoning behind the discovery and its generalizations. We describe the motivation behind the discovery of this method by setting up a symbolic experiment, observing patterns, and conjecturing a formula. We prove convergence for the square root case and establish optimal parameter conditions, then explore the parameter space to provide rigorous computational evidence.

After developing the binomial biroot, we formulate the Gaussian Biroot conjecture. Computational analysis reveals remarkably strong convergence properties. When we analyze the coefficients of Newton's method and Padé approximations through the lens of Gaussian distributions, we find compelling evidence that may explain why the biroot method works so well.

Finally, we explore other Pascal's triangle-like DAG structures to show that they also serve as approximation foundations for nth roots using a formula similar to the Binomial Biroot. We also investigate invariance properties, such as how diagonals can be used in place of rows. This structural invariance suggests an underlying mathematical structure behind nth root approximation.

\subsection{The Binomial Biroot}
\begin{figure}[H]
    \centering
    \includegraphics[width=0.75\linewidth]{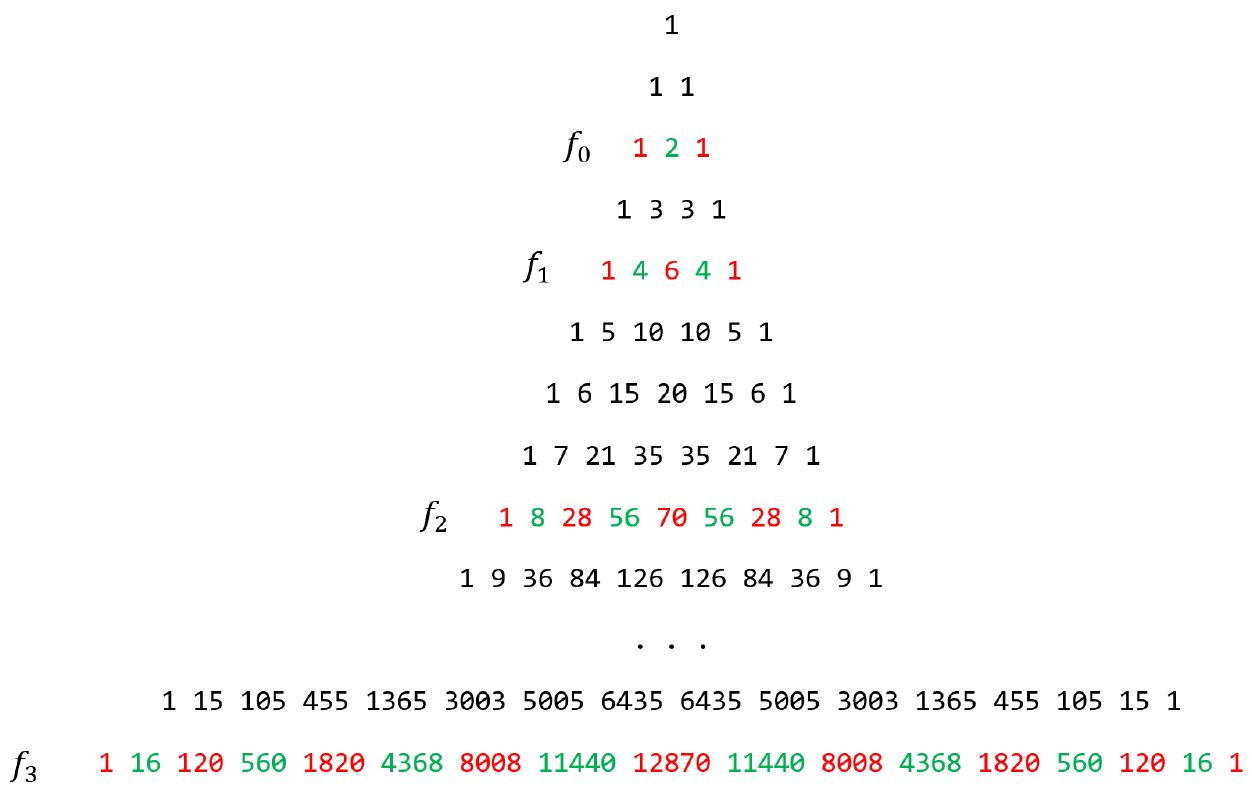}
    \caption{Pascal's triangle with highlighted rows. It should be noted that each recursive step $k$ takes us to the $2^{k+1}$th row of Pascal's triangle.}
    \label{fig:pascals_triangle_rows_from_newtons_method}
\end{figure}

We start by asking the following question: what happens if we recursively evaluate our form of Newton's method for square roots symbolically? To find out, we take the formula for Newton's method applied to nth roots and set $n$ to $2$.
\[x_{k+1} = \frac{(n-1)x_k^n + a}{nx^{n-1}_k} =\frac{x_k^2 + a}{2x_k}\]
Now we set $a$ to $x$ and $x_0$ to $1$. This gives us our first function:
\[f_0(x) = \frac{1 + x}{2}\]
We then feed this function back into into the formula: set $a$ to $x$, $x_k$ to $f_{k-1}$, simplify, and repeat.
\begin{align*}
f_1(x) &= \frac{f_0^2 + x}{2f_0} = \frac{1+6x+1x^2}{4+4x} \\
f_2(x) &= \frac{f_1^2 + x}{2f_1} = \frac{1+28x+70x^2+28x^3+1x^4}{8+56x+56x^2+8x^3} \\
f_3(x) &= \frac{f_2^2 + x}{2f_2} = \frac{1+120x+1820x^2+8008x^3+12870x^4+8008x^5+1820x^6+120x^7+ 1x^8}{16+560x+4368x^2+11440x^3+11440x^4+4368x^5+560x^6+16x^7}
\end{align*}

By now, we should notice a pattern in both the top and bottom coefficients: they appear to be binomial coefficients. More specifically, if we consider Figure~\ref{fig:pascals_triangle_rows_from_newtons_method}, we see that the coefficients are coming from rows in Pascal's triangle that alternate between the numerator and denominator.

Furthermore, when doing a similar analysis on the Padé approximations of the square root expanded at 1,
\begin{align*}
P^{1}_{1/1} &=  \frac{1 + 3 x}{3 + 1 x} \\
P^{1}_{2/1} &=  \frac{1 + 6 x + 1 x^2}{4 + 4 x} \\
P^{1}_{2/2} &=  \frac{1 + 10 x + 5 x^2}{5 + 10 x + 1 x^2} \\
P^{1}_{3/2} &=  \frac{1 + 15 x + 15 x^2 + 1 x^3}{6 + 20 x + 6 x^2} \\
P^{1}_{3/3} &=  \frac{1 + 21 x + 35 x^2 + 7 x^3}{7 + 35 x + 21 x^2 + 1 x^3} \\
P^{1}_{4/3} &=  \frac{1 + 28 x + 70 x^2 + 28 x^3 + 1 x^4}{8 + 56 x + 56 x^2 + 8 x^3} \\
\end{align*}
it appears to hold the same, but more detailed, pattern (every PT row is represented). After recognizing these patterns, we can conjecture a closed-form formula to generate these functions. It appears that by computing the binomial coefficients, we are iterating through the recursive steps of Newton’s method and/or creating Padé approximations of the square root.

\begin{remark}
Note that $f_2(x)$ from Newton's method is identical to $P^{1}_{4/3}$, demonstrating that the 3rd Newton iteration produces the [4/3] Padé approximant. It is interesting that Newton's method appears to be generating Padé approximations under the hood, hence conjecture~\ref{cjt:pade_equality}. This is further evidenced in the $c$ expansion subsection. The appearance of binomial coefficients in both Newton iterations and Padé approximants suggests these methods access the same underlying mathematical structure for rational approximation of square root functions.
\end{remark}

\begin{conjecture}\label{cjt:pade_equality}
\[P_{[p/q]}^{(c)}(x) \sim \sqrt{x}  = \beta^2_m{(x,\sqrt{c})}\]
where
\[q \in {\{p,p-1\}},  m= p+q+1\]
\end{conjecture}

\subsubsection{Closed-form formula}
Using the "n choose k" formula
\[
\binom{n}{k} = \frac{n!}{k!(n-k)!}
\]
we construct our first version of the biroot formula.
\[
\beta^2_m{(x)} = \sum_{k=0}^{m}{x^k\binom{2m}{2k}} \bigg / \sum_{k=0}^{m-1}{x^k\binom{2m}{2k+1}}
\]

\subsubsection{The $c$ Expansion Parameter (a.k.a Centering Parameter)}
To generalize our initial setup of Newton's method, instead of starting with an initial guess of 1, we start with an initial guess of c.
\begin{align*}
f_0(x) &= \frac{c^2 + x}{2c} \\
f_1(x) &= \frac{f_0^2 + x}{2f_0} = \frac{1c^4 + 6c^2x + 1x^2}{4c^3 + 4cx} \\
f_2(x) &= \frac{f_1^2 + x}{2f_1} = \frac{1c^8 + 28c^6x + 70c^4x^2 + 28c^2x^3 + 1x^4}{8c^7 + 56c^5x + 56c^3x^2 + 8cx^3}
\end{align*}
It is also important to note that the Padé expansion exhibits the same expansion pattern,
\begin{align*}
P^{(c^2)}_{1/1} &=  \frac{ 1 c^{3} + 3c x}{3 c^{2} + 1 x} \\
P^{(c^2)}_{2/1} &=  \frac{1c^{4} + 6 c^{2} x + 1x^{2}}{4c^{3} + 4c x} \\
P^{(c^2)}_{2/2} &=  \frac{1c^{5} + 10 c^{3} x + 5c x^{2}}{5 c^{4} + 10 c^{2} x + 1x^{2}}
\end{align*}
From the above, we see that as the index of $x$ increases, the index of $c$ appears to decrease in steps of two. Again, we see the index of $c$ following a similar pattern to what we have already seen: alternating even and odd indexes between the numerator and denominator. Incorporating this pattern into our closed-form formula yields
\[
\beta^2_m{(x,c)} = \sum_{k=0}^{m}{x^k c^{m-2k} \binom{2m}{2k}} \bigg / \sum_{k=0}^{m-1}{x^k c^{m-2k-1}  \binom{2m}{2k+1}}
\]

Here, it is important to note that the $c$ parameter appears to let us expand (or center) the approximation at $c^2$. This is because any initial guess for Newton's method will require no iteration if $c=\sqrt{x}$. Because $\sqrt{c^2}=c$, we know that for any initial guess of $c$, the approximation should be exact at $c^2$. Thus we may establish a property:
\[\beta^2{(c^2,c)} = c\]
This property may prove important for proofs later on, as will be seen in later sections.

While our current formula is nice in the sense that $m$ lets us control the \textit{degree} of the expression, we may want to further generalize it so that \textit{any chosen row} of Pascal's triangle can be used rather than only the $2m$th rows. To do this, we simply substitute $2m$ for $m$. If we do this, however, we must remember to change the upper bound from $m$ to $\lceil m/2 \rceil$ (although this is not strictly necessary as the terms are zero when \(2k>m\)). This yields
\[
\beta^2_m{(x,c)} = \sum_{k=0}^{\lceil \frac{m}{2} \rceil}{x^k c^{m-2k} \binom{m}{2k}} \bigg / \sum_{k=0}^{\lceil \frac{m}{2} \rceil - 1}{x^k c^{m-2k-1}  \binom{m}{2k+1}}
\]

\begin{remark}
While we don't prove the above formula through Newton's method, we do prove that it converges to the square root function using combinatorics and the limit. In this way, we essentially work backwards to confirm that it does indeed converge as expected.
\end{remark}

\begin{lemma}[Binomial Identity for Even and Odd Terms]\label{lem:binomial_even_odd_identity}
For any real number $t$ and positive integer $m$:
\begin{align}
(1 + t)^m + (1 - t)^m &= 2 \sum_{j=0}^{\lceil m/2 \rceil} \binom{m}{2j} t^{2j} \\
(1 + t)^m - (1 - t)^m &= 2t \sum_{j=0}^{\lceil (m-1)/2 \rceil} \binom{m}{2j+1} t^{2j}
\end{align}
\end{lemma}
\begin{proof}
By the binomial theorem:
\begin{align*}
(1 + t)^m &= \sum_{j=0}^m \binom{m}{j} t^j \\
(1 - t)^m &= \sum_{j=0}^m \binom{m}{j} (-1)^j t^j
\end{align*}

\noindent Adding these expansions:
\[(1 + t)^m + (1 - t)^m = \sum_{j=0}^m \binom{m}{j} t^j [1 + (-1)^j]\]
Since $[1 + (-1)^j] = 2$ when $j$ is even and $0$ when $j$ is odd, we obtain the \textit{first identity}.
\newline

\noindent Subtracting the expansions:
\[(1 + t)^m - (1 - t)^m = \sum_{j=0}^m \binom{m}{j} t^j [1 - (-1)^j]\]
Since $[1 - (-1)^j] = 2$ when $j$ is odd and $0$ when $j$ is even, factoring out $t$ from the odd powers gives the \textit{second identity}.
\end{proof}

\begin{theorem}[Square Root Case]\label{thm:square_root}
For any $x > 0$ and $c > 0$,
$$\lim_{m \to \infty} \beta^2_m(x, c) = \sqrt{x}$$
where $\beta^n_m(x,c)$ is defined as above.
\end{theorem}
\begin{proof}
\textbf{Step 1:}
Factoring out $c$ in the numerator of the biroot and substituting $u = x/c^2$ yields:
\[c^m \cdot \sum_{k=0}^{\lceil m/2 \rceil} u^k \binom{m}{2k}\]
By Lemma~\ref{lem:binomial_even_odd_identity} with $t = \sqrt{u}$:
\[c^m \cdot \sum_{k=0}^{\lceil m/2 \rceil} u^k \binom{m}{2k} = c^m \cdot \frac{(1 + \sqrt{u})^m + (1 - \sqrt{u})^m}{2}\]

\textbf{Step 2:}
Similarly, factoring out $c$ in the denominator and substituting $u = x/c^2$ yields:
\[c^{m-1} \cdot \sum_{k=0}^{\lceil m/2 \rceil - 1} u^k \binom{m}{2k+1}\]
By Lemma~\ref{lem:binomial_even_odd_identity} with $t = \sqrt{u}$:
\[c^{m-1} \cdot \sum_{k=0}^{\lceil m/2 \rceil - 1} u^k \binom{m}{2k+1} = c^{m-1} \cdot \frac{(1 + \sqrt{u})^m - (1 - \sqrt{u})^m}{2\sqrt{u}}\]

\textbf{Step 3:}
Constructing the ratio from Steps 1 and 2:
\[\beta^2_m(x, c) = c \sqrt{u} \cdot \frac{(1 + \sqrt{u})^m + (1 - \sqrt{u})^m}{(1 + \sqrt{u})^m - (1 - \sqrt{u})^m}\]

Taking the limit as $m \to \infty$: Let $a = 1 + \sqrt{u}$ and $b = 1 - \sqrt{u}$. Since $u > 0$, we have $a > 1$ and $|b/a| < 1$. Factoring out $a^m$:
\[\frac{a^m + b^m}{a^m - b^m} = \frac{1 + (b/a)^m}{1 - (b/a)^m} \to \frac{1 + 0}{1 - 0} = 1 \text{ as } m \to \infty\]

Therefore:
\[\lim_{m \to \infty} \beta^2_m(x, c) = c \sqrt{u} \cdot 1 = c\sqrt{\frac{x}{c^2}} = \sqrt{x}\]
\end{proof}

\subsubsection{Generalizing to the nth root}
The above closed-form formula works only for square roots because we discovered it using Newton's method applied to the square root. But what about cube roots, or nth roots in general? It may seem like the same process could be used, but upon further investigation, we find that this is not as straightforward as it may seem. For instance, when applying Newton's method to the cube root, the iterations are
\begin{align*}
f_0(x) &= \frac{2 + x}{3} \\
f_1(x) &= \frac{2f_0^3 + x}{3f_0^2} =  \frac{16 + 51 x + 12 x^2 + 2 x^3}{36 + 36 x + 9 x^2} \\
...
\end{align*}
The same is also true of Padé approximations. For instance the Padé expansions of the cube root yields
\begin{align*}
P^{1}_{1/1} &=  \frac{1 + 2 x}{2 + 1x} \\
P^{1}_{2/1} &=  \frac{5 + 20 x + 2 x^2}{12 + 15 x} \\
P^{1}_{2/2} &=  \frac{5 + 35 x + 14 x^2}{14 + 35 x + 5 x^2} \\
P^{1}_{3/2} &=  \frac{20 + 210 x + 168 x^2 + 7 x^3}{63 + 252 x + 90 x^2} \\
...
\end{align*}
While it is interesting that the coefficients are symmetrical for even expansions, our analysis of the patterns has not yet resulted in a rule or solid pattern that can be used to generate higher degree expansions. However, both the Newton and Padé expansions for nth roots do yield noteworthy results that provide insight as to why our generalized nth root biroot formula works so well. Further analysis of Newton and Padé coefficients reveals they cluster in pairs of two around a Gaussian distribution.

Despite this complexity, computational experimentation reveals a surprising formula: while direct symbolic iteration fails to produce simple formulas, sampling every nth pair of coefficients from a row in Pascal's triangle yields remarkably accurate $n$th root approximations. Essentially, we simply substitute $2$ for $n$ (where $n$ is the index of the root) in our current biroot formula and we get the following generalized formula:
\[
\beta^n_m{(x,c)} = \sum_{k=0}^{\lceil \frac{m}{n} \rceil}{x^k c^{m-nk} \binom{m}{nk}} \bigg / \sum_{k=0}^{\lceil \frac{m}{n} \rceil - 1}{x^k c^{m-nk-1} \binom{m}{nk+1}}
\]

It is worth noting that for the expansion parameter $c$, the aforementioned process of using Newton's method or Padé approximants on nth roots does indicate that the $c$ parameter's degree decreases in steps of $n$ for each term as $x$'s degree increases. This helps validate the above generalized formula by adding another line of evidence in regard to the $c$ parameter.

Thus we have the generalized form of the biroot method. We can now construct Padé-like approximations for the any nth root function. For instance, using the 12th row of Pascal's triangle to construct a cube root approximation yields
\[
\frac{1 + 220 x + 924 x^2 + 220 x^3 + 1x^4}{12 + 495 x + 792 x^2 + 66 x^3}
\]
\begin{figure}
    \centering
    \includegraphics[width=0.75\linewidth]{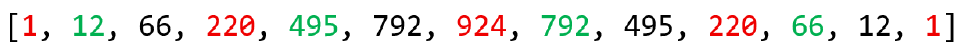}
    \caption{Red values are the coefficients of the numerator and Green values are the coefficients of denominator. Notice that coefficients are in pairs of two and are being sampled at every third position in the row.}
    \label{fig:cube_root_from_pascals_12_row}
    
\end{figure}
Figure~\ref{fig:cube_root_from_pascals_12_row} shows these coefficients coming from the 12th row of pascals triangle. This function has an error of around $0.0001$ within the interval $(0.07, 13)$.

\begin{conjecture}[Generalized Binomial Biroot]\label{cjt:generalized_binomial_biroot}
For any $x > 0$, positive integer $n$, and $c > 0$,
$$\lim_{m \to \infty} \beta^n_m(x, c) = \sqrt[n]{x}$$
where $\beta^n_m(x,c)$ is defined as above.
\end{conjecture}

We do not attempt a general proof of Conjecture~\ref{cjt:generalized_binomial_biroot} in this work as that is beyond the scope of this paper. Instead, we focus on establishing the biroot method and its apparent fundamental properties. While the square root case demonstrates an elegant proof through binomial identities and standard limit analysis, the general nth root case likely requires a similar, but more sophisticated technique. A complete proof would likely employ roots of unity filters, combined with asymptotic analysis of the resulting rational expressions. Such an approach would naturally extend our binomial identity framework to handle the modular sampling pattern inherent in the general biroot construction. However, the extensive computational evidence presented throughout this paper, combined with the square root proof, provides a strong evidence-based foundation for the conjecture.

\begin{remark}
The general proof would likely proceed by: (1) applying roots of unity filters to isolate coefficients at positions $nk$, (2) expressing the Biroot's numerator and denominator as sums over nth roots of unity, (3) analyzing the dominant term of the the ratio of those sums as $m \to \infty$, and (4) showing convergence to $x^{1/n}$. This would be a similar, but generalized, technique as used in the proof for theorem \ref{thm:square_root}.
\end{remark}

\subsubsection{Computational Evidence}
\begin{figure}[H]
    \centering
    \includegraphics[width=0.75\linewidth]{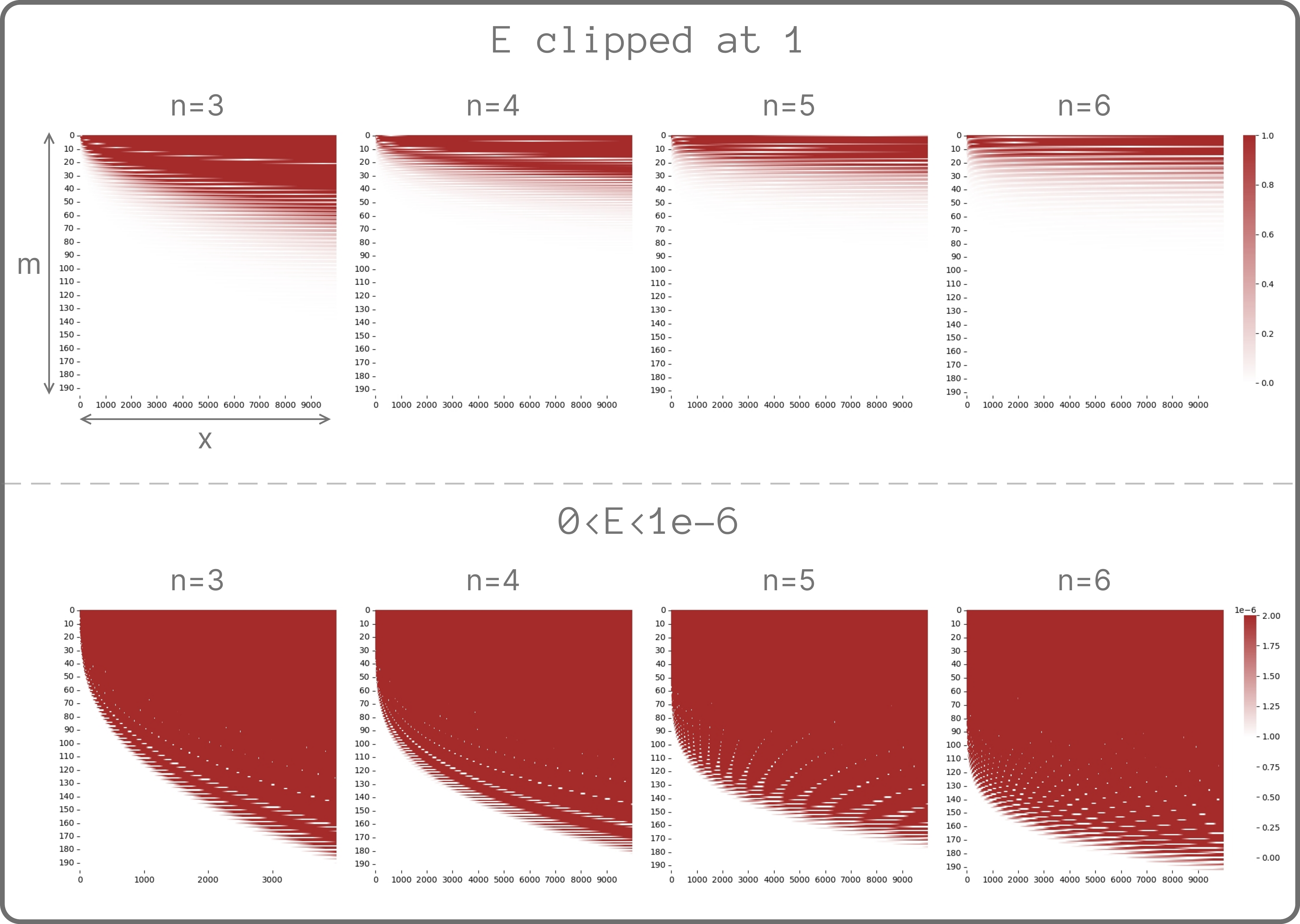}
    \caption{$c=1$: Convergence of $x$ in the interval $[0, 10^4]$ when varying $m$ over $[n+1, 200]$, and $n$ over $3, 4, 5$, and $6$.}
    \label{fig:binom_heatmaps_c_1}
\end{figure}
Using our custom computational framework (along with Seaborn and Matplotlib), we explore the parameter space of the biroot method to understand its convergence properties. This analysis serves both as evidence of convergence and as guidance for choosing optimal parameters in practical applications. We have archived the computational matrices containing approximately 24 million parameter evaluations in \verb|.npz| files to facilitate future analysis and reproducibility.

Extensive empirical analysis across the parameter space reveals promising convergence properties for the generalized biroot method. Heat map visualizations (using the error function $E=|\beta^n_m(x, c) - \sqrt[n]{x}|$) spanning $x \in [0, 10^4]$ and $m \in [n+1, 200]$ demonstrate rapid error decay as the approximation order increases, with error magnitudes decreasing by orders of magnitude as $m$ increases. The error patterns exhibit distinct regimes: initial moderate convergence for low-degree approximations, followed by dramatically accelerating convergence rates for higher-orders.

\begin{figure}[H]
    \centering
    \includegraphics[width=0.9\linewidth]{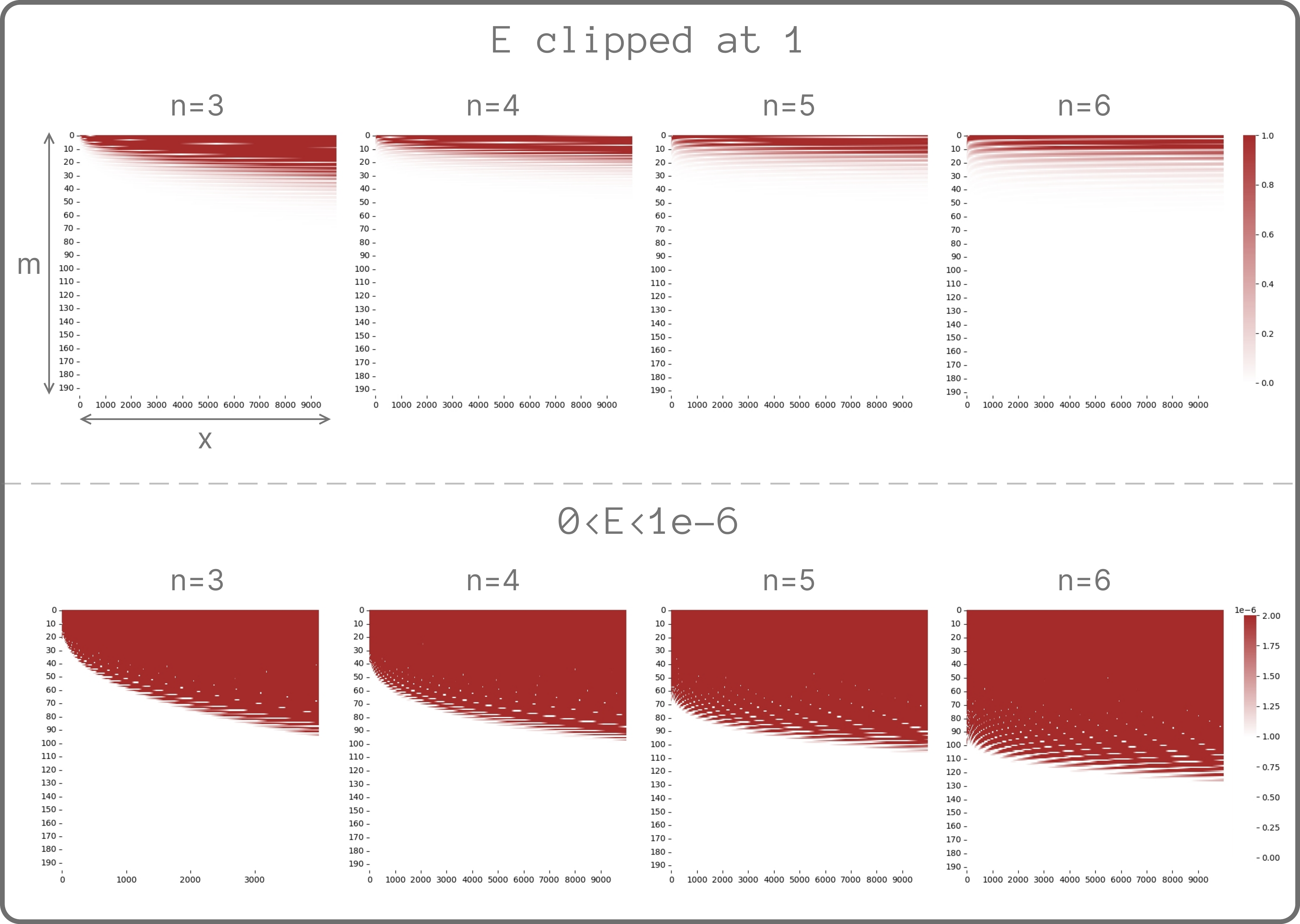}
    \caption{$c=2$}
    \label{fig:binom_heatmaps_c_2}
\end{figure}

\begin{figure}[H]
    \centering
    \includegraphics[width=0.9\linewidth]{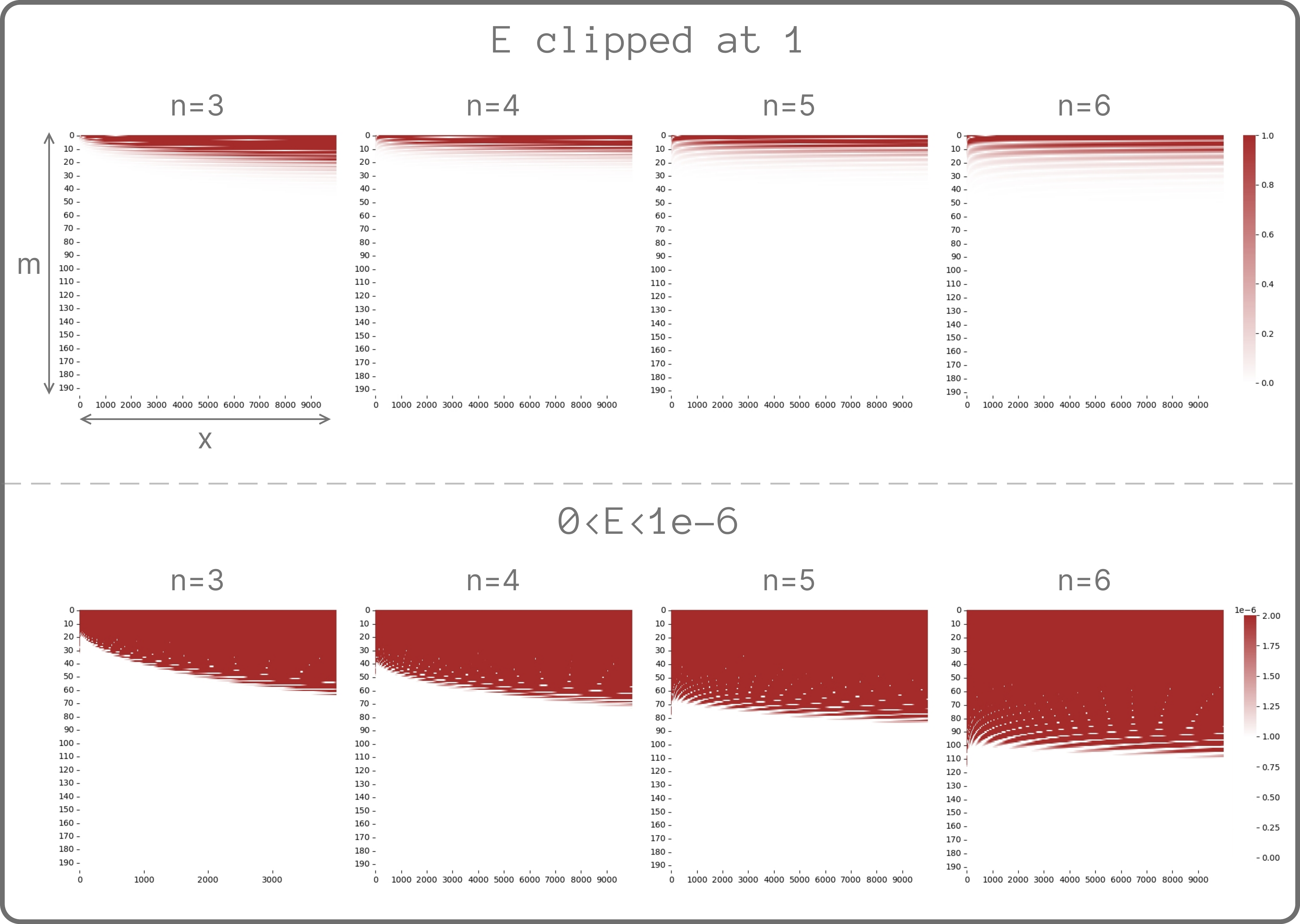}
    \caption{$c=3$}
    \label{fig:binom_heatmaps_c_3}
\end{figure}

\begin{figure}[H]
    \centering
    \includegraphics[width=0.9\linewidth]{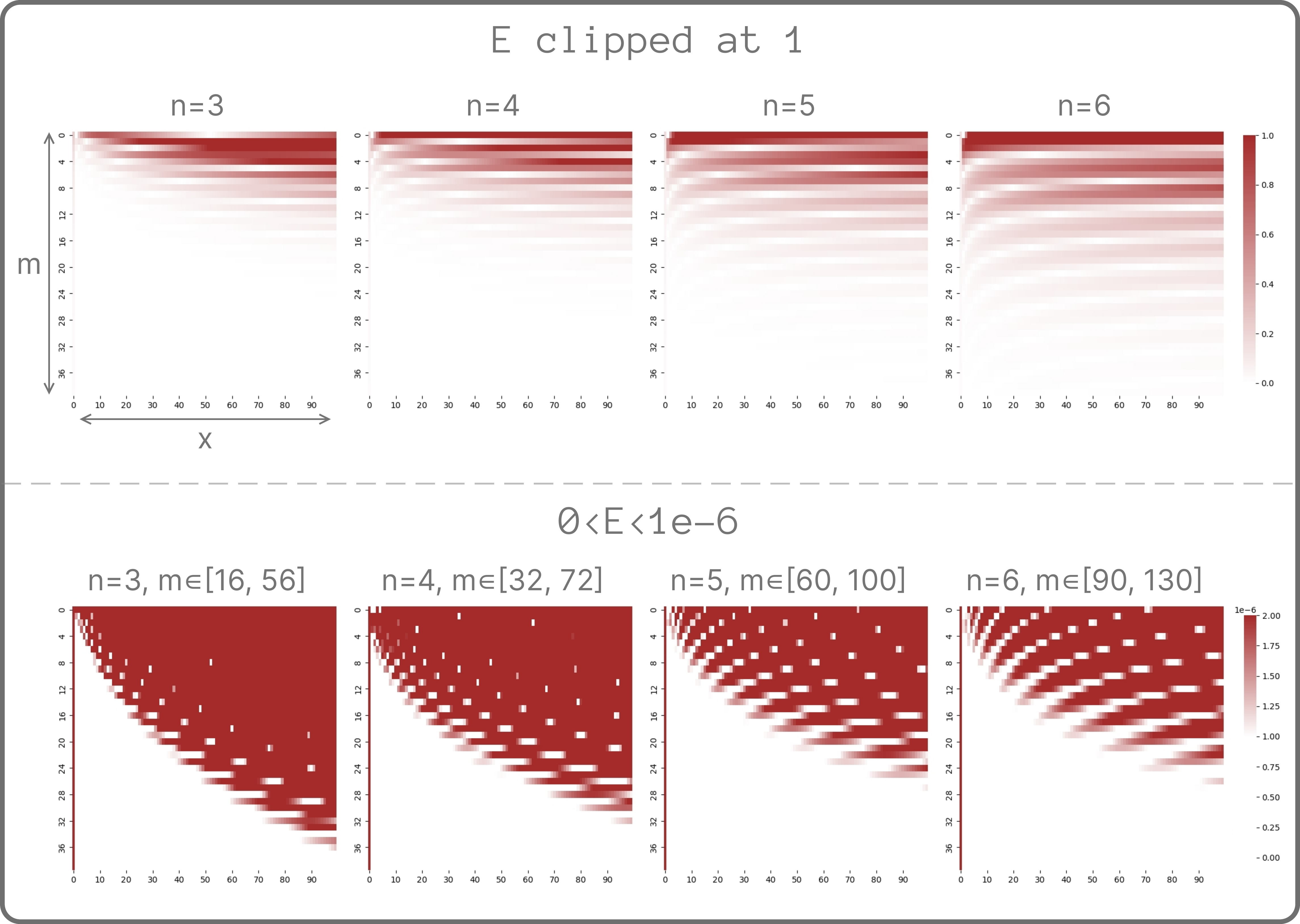}
    \caption{Zoomed in with $c=1$: Convergence of $x$ in the interval $[0, 100]$ when varying $m$ over $[n+1, 40]$ (unless otherwise specified), and $n$ over $3, 4, 5$, and $6$.}
    \label{fig:binom_heatmaps_pv_c_1}
\end{figure}

The influence of the expansion parameter $c$ proves particularly significant, with larger values of $c$ yielding substantially improved convergence rates across all tested configurations. Heatmap \cref{fig:binom_heatmaps_c_2,fig:binom_heatmaps_c_3,fig:binom_heatmaps_pv_c_1} for $c=2, 3$ and a zoomed in version are available in the appendices.

It should be noted that the included figures show error tolerance only for $0$ to $1$ or $0$ to $10^{-6}$ (high precision tolerance). When using higher error tolerance thresholds, one can obtain significantly wider intervals of convergence. This may be useful, for instance, to obtain seed values that can then be refined using an iteration of Newton's method or other root-finding techniques. This consideration is also relevant to the implementation of the recursive biroot algorithm (discussed later).

Preliminary log-log regression analysis over the parameter ranges $m \in [9, 79]$ and log-spaced interval $x \in [0.1, 2048]$ suggests convergence rates potentially following power-law $O(m^{-\alpha})$ behavior. Estimated exponents $\alpha$ range from approximately $12$ to $17$ for the cube root case ($n=3, c=1$), and $8$ to $13$ for the fourth root case ($n=4, c=1$). While these preliminary estimates indicate exceptionally rapid convergence compared to classical rational approximation methods, the substantial variability across the parameter space indicates that convergence behavior depends systematically on the relationship between $x$, $c$, and $n$, warranting more extensive theoretical analysis.

The computational evidence strongly supports the validity of Conjecture~\ref{cjt:generalized_binomial_biroot}, demonstrating that the biroot method consistently converges to the nth root function across diverse parameter ranges. The systematic error patterns observed in the heat map analysis reveal fundamental structural properties that merit extensive theoretical investigation.

\subsection{Optimal Condition Theorem}
A natural question regarding our generalized biroot formula concerns the expansion property: for which values of $m$ and $n$ does $\beta^n_m(c^n,c) = c$ hold exactly? To investigate this, we evaluate $\beta^n_m(c^n,c)$:
\begin{align*}
\beta^n_m{(c^n,c)} &= \sum_{k=0}^{\lceil \frac{m}{n} \rceil}{c^{nk} c^{m-nk} \binom{m}{nk}} \bigg / \sum_{k=0}^{\lceil \frac{m}{n} \rceil - 1}{c^{nk} c^{m-nk-1}  \binom{m}{nk+1}} \\ 
&= \sum_{k=0}^{\lceil \frac{m}{n} \rceil}{c^m \binom{m}{nk}} \bigg / \sum_{k=0}^{\lceil \frac{m}{n} \rceil - 1}{c^{m-1}  \binom{m}{nk+1}} \\
&= c \cdot \left( \sum_{k=0}^{\lceil \frac{m}{n} \rceil}{\binom{m}{nk}} \bigg / \sum_{k=0}^{\lceil \frac{m}{n} \rceil - 1}{\binom{m}{nk+1}} \right)
\end{align*}
We then use this result to define the coefficient ratio function:
\[
\alpha(m, n) = \frac{\sum_{k=0}^{\lceil \frac{m}{n} \rceil}{\binom{m}{nk}}}{\sum_{k=0}^{\lceil \frac{m}{n} \rceil - 1}{\binom{m}{nk+1}}}
\]

Therefore, the fixed-point property $\beta^n_m(c^n,c) = c$ holds if and only if $\alpha(m,n) = 1$. To ascertain when the above expression appears to equal exactly $1$, we use our computational framework to run an analysis. This way, we may be able to modify the expression so that it always equals 1, and then prove it. Some of the results from our analysis are shown in table 1.

\begin{table}[H]
    \centering
    \small
    \begin{tabular}{c|c|c|c|c|c|c|c|c|c|c|c|c}
         $n$ & \multicolumn{12}{c}{$m$ values where $\alpha(m,n) = 1$} \\
         \hline
         3 & 4 & 7 & 10 & 13 & 16 & 19 & 22 & 25 & 28 & 31 & 34 & 37 \\
         4 & 5 & 9 & 13 & 17 & 21 & 25 & 29 & 33 & 37 & 41 & 45 & 49 \\
         5 & 6 & 11 & 16 & 21 & 26 & 31 & 36 & 41 & 46 & 51 & 56 & 61 \\
         6 & 7 & 13 & 19 & 25 & 31 & 37 & 43 & 49 & 55 & 61 & 67 & 73 \\
         \hline
         \multicolumn{13}{c}{Pattern: $m \equiv 1 \pmod{n}$} \\
    \end{tabular}
    \caption{Optimal parameter combinations for exact fixed-point property}
    \label{tab:optimal_pattern}
\end{table}

Here we see that when $m \equiv 1 \pmod{n}$, then $\alpha(m,n)$ is exactly equal to $1$. Thus, we modify our coefficient ratio function so that this is always true. We do this by substituting $m$ for $mn+1$:
\[
\alpha^{'}(m, n) = \frac{\sum_{k=0}^{m+1}{\binom{mn+1}{nk}}}{\sum_{k=0}^{m}{\binom{mn+1}{nk+1}}}
\]
We use $mn+1$ because it enforces the condition $m \equiv 1 \pmod{n}$. The upper bounds are also changed because $\lceil \frac{m}{n} \rceil$ becomes $\lceil \frac{mn+1}{n} \rceil = m+1$. We then claim that $\alpha^{'}(m, n) = 1$ for all valid $m$, $n$, which we prove below.

\begin{remark}
    Theorem~\ref{thm:fixed_point} establishes that for any root index $n \ge 2$, there are infinitely many biroot approximations $\beta^n_m{(x,c)}$ satisfying the exact fixed-point property $\beta^n_m{(c^n,c)}=c$, occurring precisely when $m \equiv 1 \pmod{n}$
\end{remark}

\begin{lemma}\label{lem:optimal_condition}
When $n \ge 2$,
\[
\sum_{k=0}^{m+1}{\binom{mn+1}{nk}} = \sum_{k=0}^{m}{\binom{mn+1}{nk+1}}
\]
\end{lemma}

\begin{proof}
Using the symmetry property $\binom{n}{k}=\binom{n}{n-k}$, the terms on the RHS become
\[\binom{mn+1}{nk+1} = \binom{mn+1}{(mn+1)-(nk+1)} = \binom{mn+1}{mn-nk}\]
So the RHS becomes
\[\sum_{k=0}^{m}{\binom{mn+1}{mn-nk}} = \sum_{k=0}^{m}{\binom{mn+1}{n(m-k)}}\]
Substituting $m-k$ for $j$ yields (this is valid because when $k$ ranges from $0$ to $m$, $j$ ranges from $m$ to $0$)
\[\sum_{k=0}^{m}{\binom{mn+1}{n(m-k)}} = \sum_{j=0}^{m}{\binom{mn+1}{nj}}\]
Now observe that the LHS of the identity can be split
\[\sum_{k=0}^{m+1}{\binom{mn+1}{nk}} = \sum_{k=0}^{m}{\binom{mn+1}{nk}} + \binom{mn+1}{n(m+1)}\]
Since we've shown that the RHS equals $\sum_{k=0}^{m}{\binom{mn+1}{nk}}$, we need
\[\sum_{k=0}^{m}{\binom{mn+1}{nk}} + \binom{mn+1}{n(m+1)} = \sum_{j=0}^{m}{\binom{mn+1}{nj}}\]
This requires that $\binom{mn+1}{n(m+1)} = 0$. Because we know that $\binom{n}{k}$ equals $0$ when $k>n$, we recognize that the identity only holds when $nm + n > mn + 1$. This inequality simplifies to $n > 1$. Thus, the the identity only holds when $n \ge 2$.
\end{proof}

\begin{theorem}\label{thm:fixed_point}
By Lemma~\ref{lem:optimal_condition},
\[\beta^n_m(c^n,c) = c\]
when $n \ge 2$, and $m \equiv 1 \pmod{n}$.
\end{theorem}

\subsection{The Gaussian Biroot}
The use of the Gaussian distribution as an approximation to the binomial distribution is well-established through the Central Limit Theorem  (more specifically the De Moivre-Laplace Theorem), a fundamental result demonstrating that binomial distributions converge to the normal distribution as the number of trials increases. This approximation is widely accepted in both theoretical and applied contexts. Hence, we ask, why not modify our current Binomial Biroot formula to make it a Gaussian-based Biroot formula? Furthermore, might this provide better convergence by being the underlying attractor of the binomial distribution?

We start by defining the normal distribution

\[G(x) =\frac{1}{\sigma \sqrt{2 \pi}} \cdot e^{- \frac{(x-\mu)^2}{2 \sigma^2}}\]

Because we are working with a binomial distribution with parameters $m$ and $p=1/2$, we must use its standard formulas for the mean $\mu$ and standard deviation $\sigma$. Applying these yields
\begin{align*}
\mu &= mp = \frac{m}{2} \\
\sigma &= \sqrt{mp(1-p)} = \frac{\sqrt{m}}{2}
\end{align*}
Because $\frac{1}{\sigma \sqrt{2 \pi}}$ is only a scalar, it is completely irrelevant when sampling coefficients from the Gaussian function. Thus, the following Gaussian function is used.
\[G(x) =e^{-\frac{2(x-\frac{m}{2})^2}{m}}\]

Now it is just a matter of sampling alternated points from the gaussian to use as coefficients. Hence, our new Gaussian Biroot formula:
\[
\beta^n_m{(x)} = \sum_{k=0}^m x^k G(nk ) \bigg / \sum_{k=0}^{m-1} x^k G(nk+1)
\]

Incorporating the $c$ expansion parameter yields
\[
\beta^n_m{(x)} = \sum_{k=0}^m x^k c^{m-nk}  G(nk ) \bigg / \sum_{k=0}^{m-1} x^k c^{m-nk-1} G(nk+1)
\]

\begin{remark}
    Although the upper bound should be $\lceil \frac{m}{n} \rceil$ rather than $m$, we find the convergence works better when left unconstrained (this requires further investigation). Both, however, do seem to work.
\end{remark}

\subsubsection{Evidence of Gaussian Coefficient Mapping}
When plotting the normalized values of a binomial distribution (row in Pascal's triangle) as points, we quickly see that as larger distributions are chosen, they converge to points on the Gaussian function defined above, where $m$ is the size of the given distribution (or row). The following figure shows such a distribution when $m=10$ compared to $m=100$:
\begin{figure}[H]
    \centering
    \includegraphics[width=0.75\linewidth]{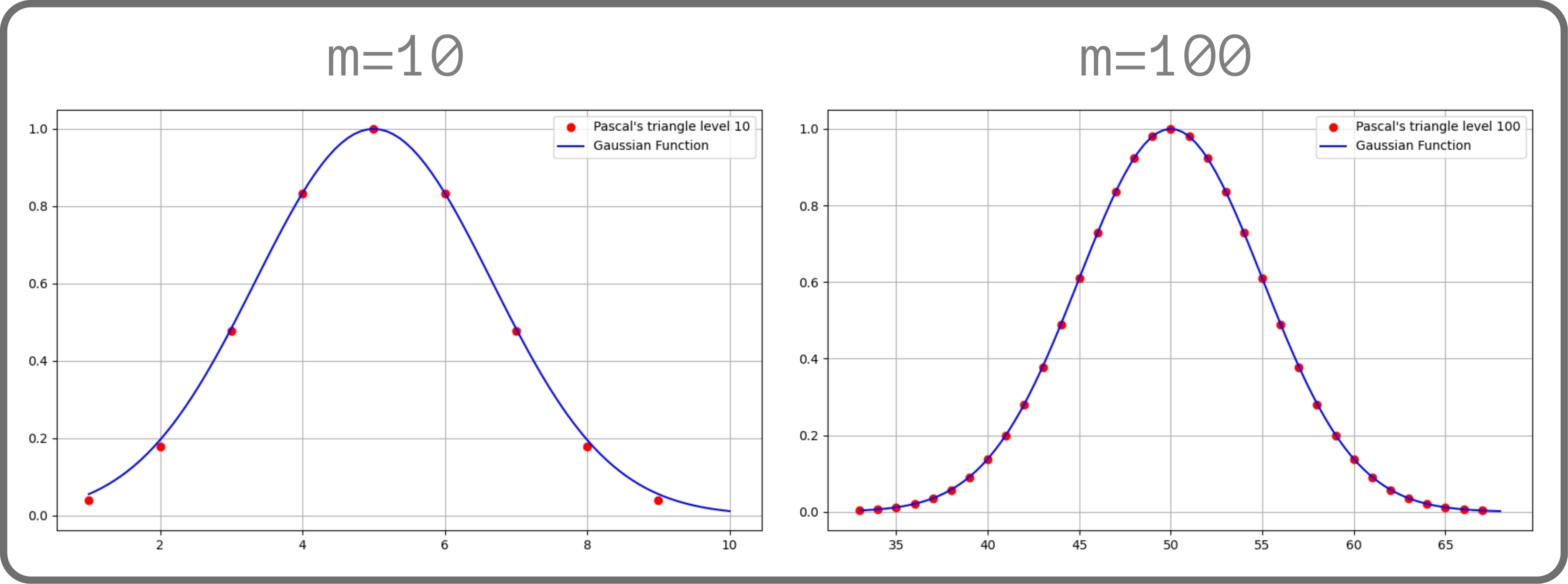}
\end{figure}

Furthermore, when plotting the coefficient distributions from Newton's method applied to nth roots, we also see the points converge to points on the Gaussian function. The following figure shows such a distribution comparing the \textit{second} symbolic iteration to the \textit{third} symbolic iteration for both the \textit{cube root} and \textit{fourth root}:
\begin{figure}[H]
    \centering
    \includegraphics[width=0.75\linewidth]{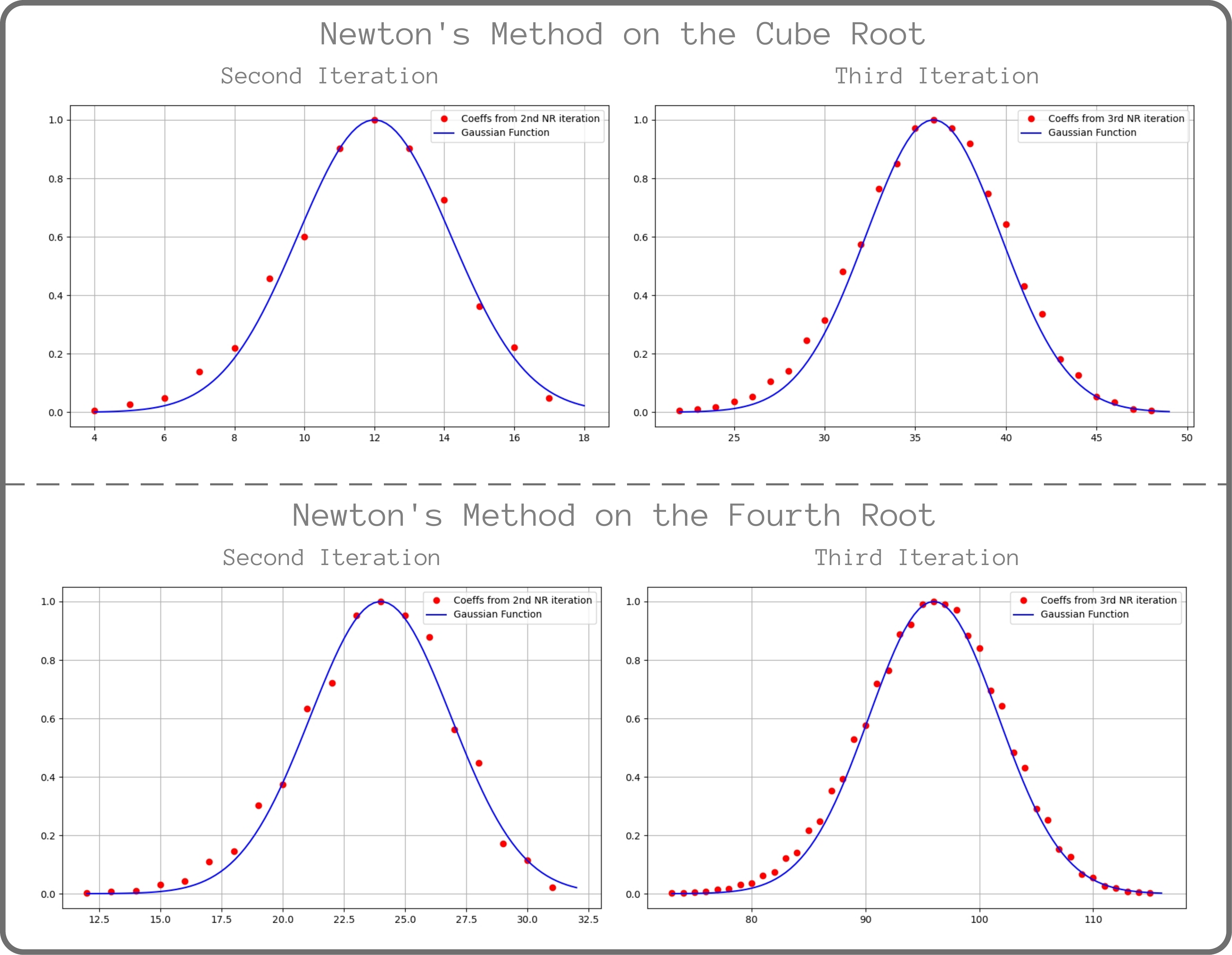}
\end{figure}

Likewise, when plotting the coefficient distributions of Padé approximants for nth roots, we see the same pattern: points converge to points on the Gaussian function:
\begin{figure}[H]
    \centering
    \includegraphics[width=0.75\linewidth]{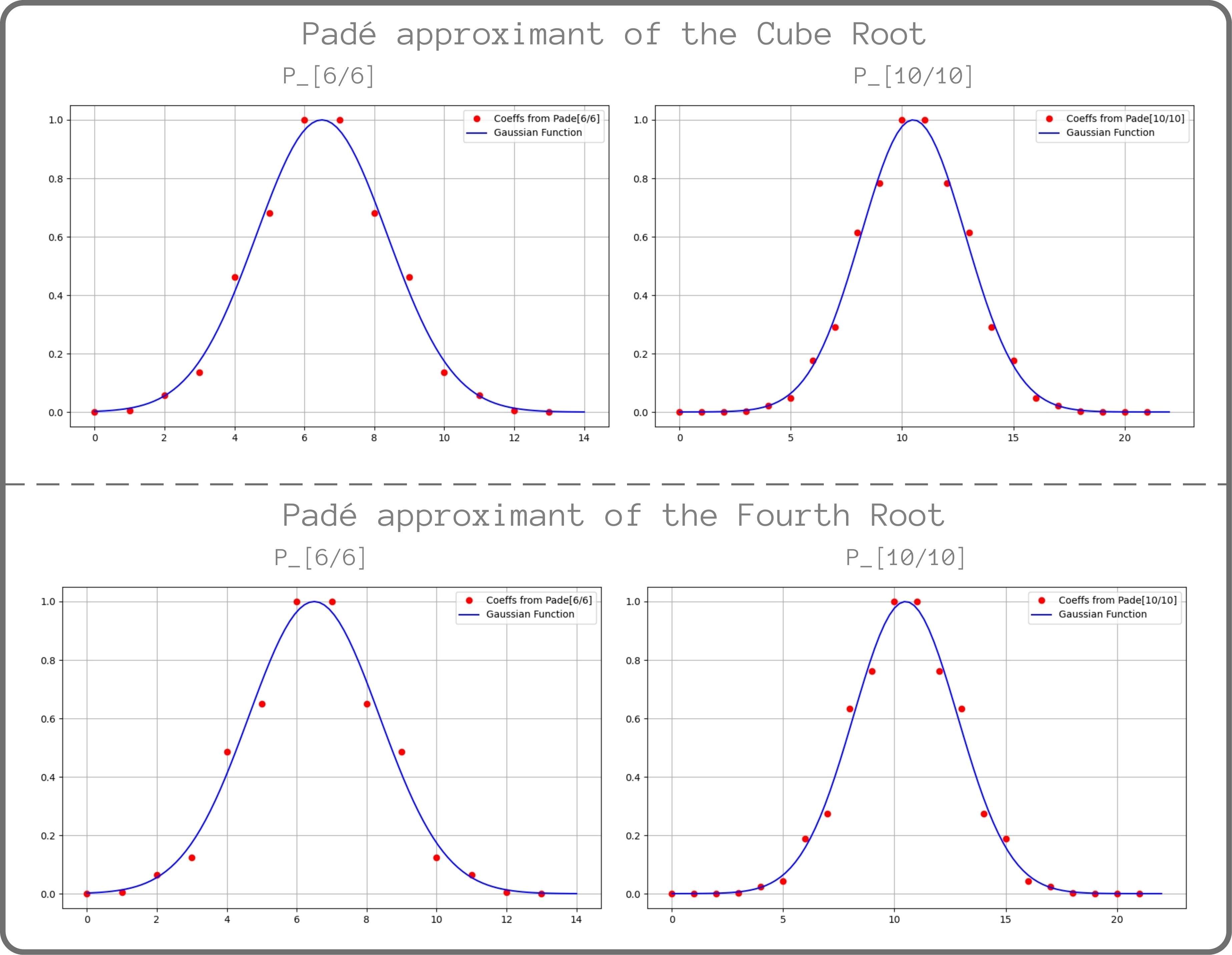}
\end{figure}

\noindent Both of the above figures demonstrate another interesting property for nth roots: the coefficients of the numerator and denominator in both Newton’s method and Padé approximants appear to follow two distinct Gaussian distributions, with the Gaussian curves separated by a small delta. This, of course, warrants further theoretical investigation.

\subsubsection{Computational Evidence}
The Gaussian Biroot exhibits remarkable convergence properties that significantly exceed those of the Binomial Biroot. Preliminary analysis suggests exponentially faster convergence rates in many parameter regimes, providing strong evidence for the Gaussian distribution as an underlying attracting structure for rational nth root approximations. However, we emphasize that this represents initial computational exploration, and more rigorous theoretical and numerical analysis is needed to fully establish these claims.

\begin{figure}[H]
    \centering
    \includegraphics[width=0.75\linewidth]{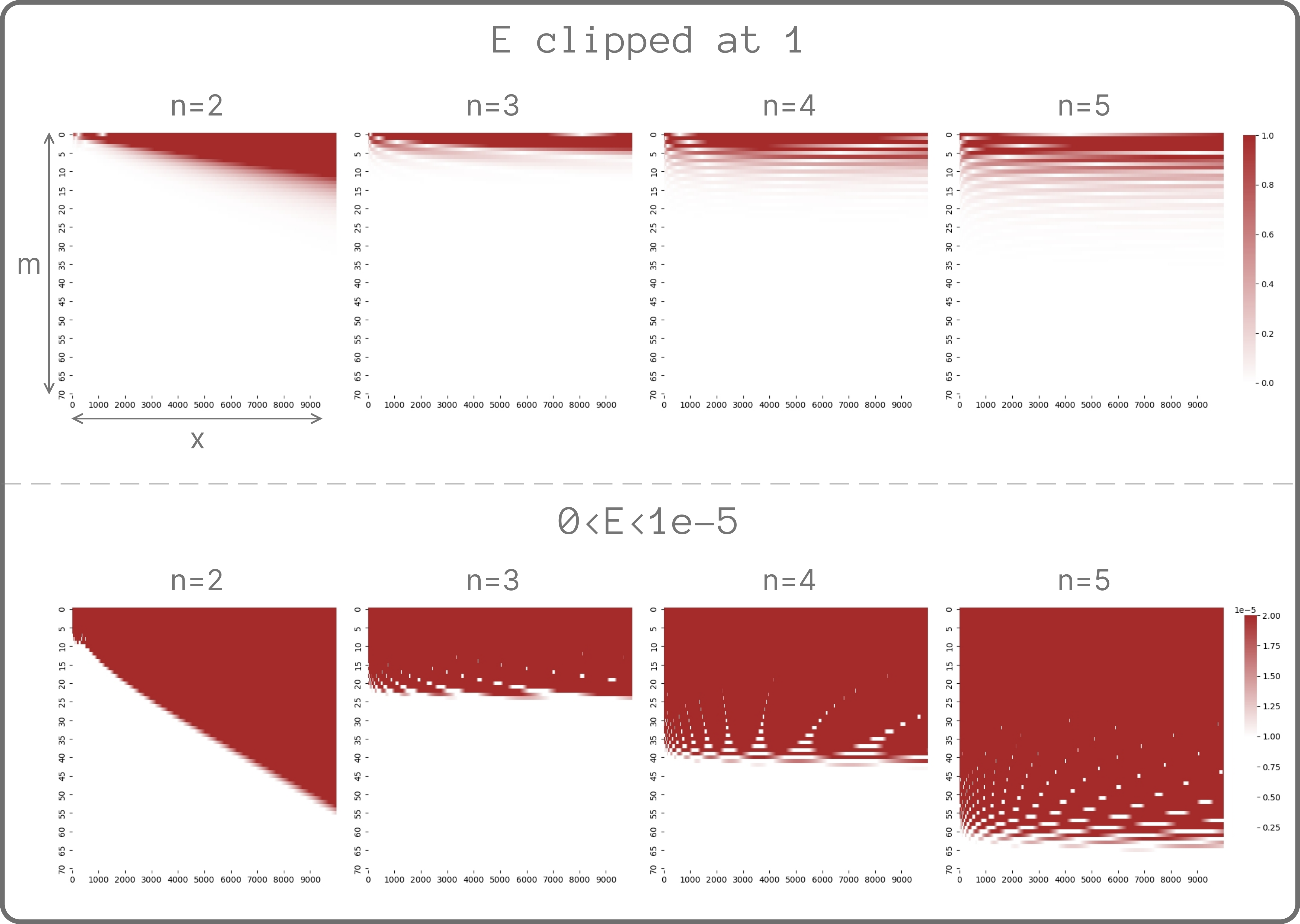}
    \caption{Gaussian Biroot convergence analysis with $c=1$. Note the significantly faster convergence compared to Binomial Biroot results, achieved with lower maximum order ($m \leq 70$ versus $m \leq 200$).}
    \label{fig:gaussian_heatmaps_c_1}
\end{figure}

The heatmaps in \cref{fig:gaussian_heatmaps_c_1} demonstrate these convergence properties using the error function $E=|\beta^n_m(x, c) - \sqrt[n]{x}|$. The approximation order $m$ varies over $[0, 70]$ while $x$ spans $[0, 10^4]$.

These results are particularly notable because they achieve superior convergence rates with substantially lower degrees compared to the Binomial Biroot. While our Binomial analysis required $m \in [0, 200]$, the Gaussian variant achieves comparable or better accuracy with $m \in [0, 70]$. Due to computational limitations in our framework's C backend (specifically overflow errors at higher values), we restricted the domain to $x \in [0, 10^4]$. However, preliminary tests suggest the convergence advantages extend to much larger domains.

To quantify performance over extended ranges, Table~\ref{tab:mean_gaussian_convergence_error} presents statistical analysis for the cube root case ($n=3, c=1$) over the interval $x \in [0, 10^8]$. For each value of $m$, we sampled $10^7$ uniformly distributed points and computed the mean absolute error $\mu$ and standard deviation $\sigma$.

The results in Table~\ref{tab:mean_gaussian_convergence_error} demonstrate the practical effectiveness of the Gaussian Biroot approach. With degree of just $m=12$, the method achieves mean absolute error below 1.1 across the interval $[0, 10^8]$ with relatively low variance ($\sigma \approx 0.53$). By degree $m=18$, the mean error drops to approximately 0.035 with standard deviation 0.021, indicating both high accuracy and consistency across the domain.

\begin{table}[H]
    \centering
    \small
    \begin{tabular}{c|ccc}
         $m$ & Mean Error ($\mu$) & Std. Dev. ($\sigma$)  &Max Error\\
         \hline
         12 & 1.0929 & 0.5279  &1.8115\\
         13 & 0.5312 & 0.2905  &1.0543\\
         14 & 0.4105 & 0.2242  &0.7019\\
         15 & 0.1917 & 0.1101  &0.3533\\
         16 & 0.1106 & 0.0616  &0.2297\\
         17 & 0.0783 & 0.0416  &0.1522\\
         18 & 0.0352 & 0.0212  &0.1512\\
 19& 0.0226& 0.0128&0.1503\\
 20& 0.0147& 0.0076&0.1495\\
 21& 0.0067& 0.0039&0.1488\\
 22& 0.0045& 0.0025&0.1482\\
    \end{tabular}
    \caption{Statistical convergence analysis for Gaussian Biroot cube root approximation over $x \in [0, 10^8]$ with $10^7$ sample points per order $m$.}
    \label{tab:mean_gaussian_convergence_error}
\end{table}

Similar performance characteristics have been observed for fourth and fifth roots, though we omit detailed results for brevity. The systematic error reduction as $m$ increases, combined with the low variance across large input ranges, supports our conjecture regarding the Gaussian distribution's fundamental role in rational root approximation.

Even when high precision tolerance is not required, the Gaussian Biroot's ability to generate accurate approximations with low-degree rational functions has potential practical value. These expressions could serve as high-quality seed values for iterative refinement using traditional numerical methods, potentially reducing the computational cost of high-precision root calculations.

\subsection{The DAG Biroot}
A noteworthy structural invariance property seems to persist for linearly constructed DAG structures. This stems from the following generalization of pascal-like triangular arrays. Using the DAG class from our computational framework, we can construct these structures using any initial basin of attraction, linear function, node arity (number of parent nodes to pass as arguments to the function), and step size. For instance, \verb|DAG().as_graph(basin=[1], depth=10, node_arity=2)| yields ten levels of Pascal's triangle.

A DAG starting with a random basin of \verb|[3, 1, 6, 7]|, node arity of 3, depth of 3, and linear function of $sum(x_0, x_1, x_3)$, for instance, would yield the following structure:
\begin{figure}[H]
    \centering
    \includegraphics[width=0.5\linewidth]{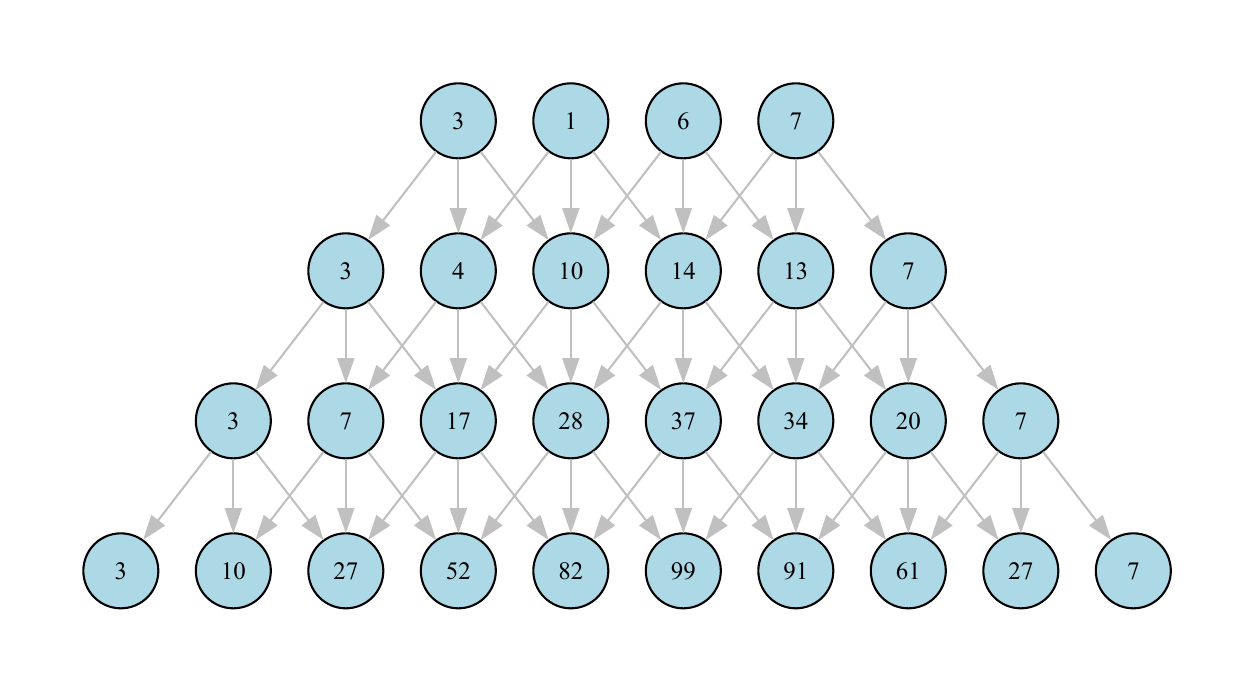}
\end{figure}

Using this generalized framework, we may ask ourselves: can any row from our generalized DAG structure be used to approximate nth roots? We do find this to be the case. However, due to the large number of configurations and test cases that would need to be explored for thorough validation, we present only preliminary computational analysis.

Consider the following structure with an arbitrary basin $[5, 2, 7, 1, 8]$ and linear function $f(a, b, c)=1a+4b+3c$:
\begin{figure}[H]
    \centering
    \includegraphics[width=0.5\linewidth]{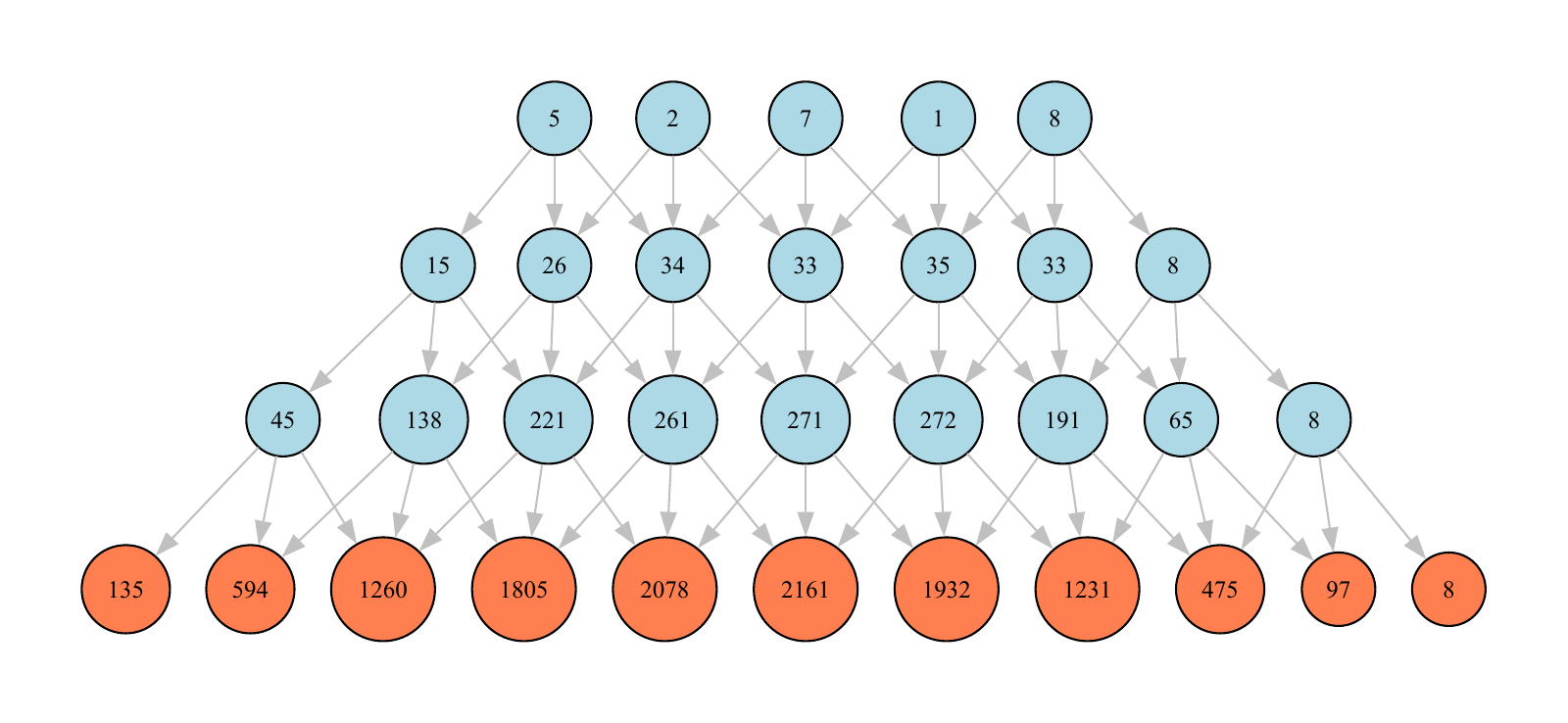}
\end{figure}
\noindent Alternating the values from the highlighted row as coefficients for a square biroot, we get
\[f(x)=\frac{135 + 1260 x + 2078 x^{2} + 1932 x^{3} + 475 x^{4} + 8 x^{5}}{594 + 1805 x + 2161 x^{2} + 1231 x^{3} + 97 x^{4}}\]
We find that $f$ has a convergence error below 0.001 in the interval $[0.25, 4.5]$.

\begin{figure}[H]
    \centering
    \includegraphics[width=0.75\linewidth]{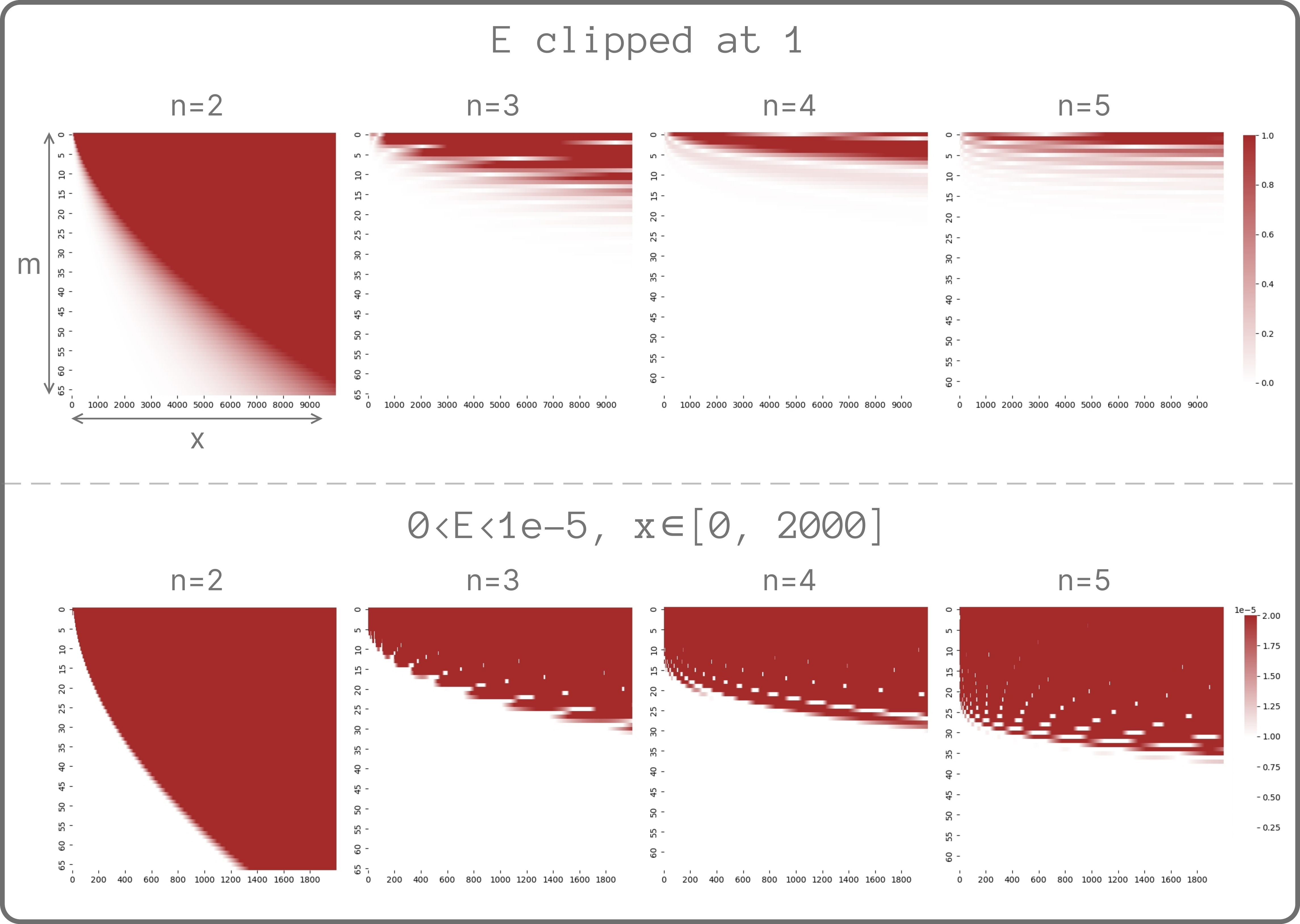}
    \caption{DAG Biroot convergence analysis with $c=1$ for the DAG with random basin $[5, 2, 7, 1, 8]$ and arbitrary node arity 3}
    \label{fig:DAG_heatmaps_c_1}
\end{figure}

The heatmaps in \cref{fig:DAG_heatmaps_c_1} demonstrate these convergence properties for the above DAG using the error function $E=|\beta^n_m(x, c) - \sqrt[n]{x}|$. The DAG's level $m$ varies over $[0, 70]$ while $x$ spans $[0, 10^4]$ ($[0, 2000]$ for error tolerance $10^{-5}$). The fact that we arbitrarily chose the basin, linear function, and node arity is strong evidence that this pattern generalizes. Nevertheless, to quantify performance over extended configurations, ~\cref{tab:mean_dag_cube_convergence_error} presents statistical analysis for the cube over the interval $x \in [0, 10^4]$. With a DAG's level $m=40$, we generate random basin, node arity, and linear function configurations. For each configuration, we sampled $10^4$ uniformly distributed points and computed the mean absolute error $\mu$ and standard deviation $\sigma$. Here we, perhaps unsurprisingly, see decent evidence of convergence. Hence, we establish the DAG-based Biroot formula:
\[
\beta^n_m{(x^n,c)} = \sum_{k=0}^{\lceil \frac{m}{n} \rceil}{x^k c^{m-nk} C(m,nk)} \bigg / \sum_{k=0}^{\lceil \frac{m}{n} \rceil - 1}{x^k c^{m-nk-1}  C(m,nk+1)}
\]
where $C(m, k)$ gets the $k$th node at level $m$ of an arbitrary linearly-constructed DAG.

\subsubsection{Gaussian Foundation}
While not rigorously justified, we hypothesize that the reason the DAG biroot works is because each level progressively approaches the Gaussian (from the previous subsection) when normalized. This would elegantly explain why it works. The following figure \ref{fig:dag_vs_gaussian} portrays how the normalized node values approach points on the Gaussian as lower levels are selected of a randomly configured DAG.

\begin{figure}[H]
    \centering
    \includegraphics[width=0.75\linewidth]{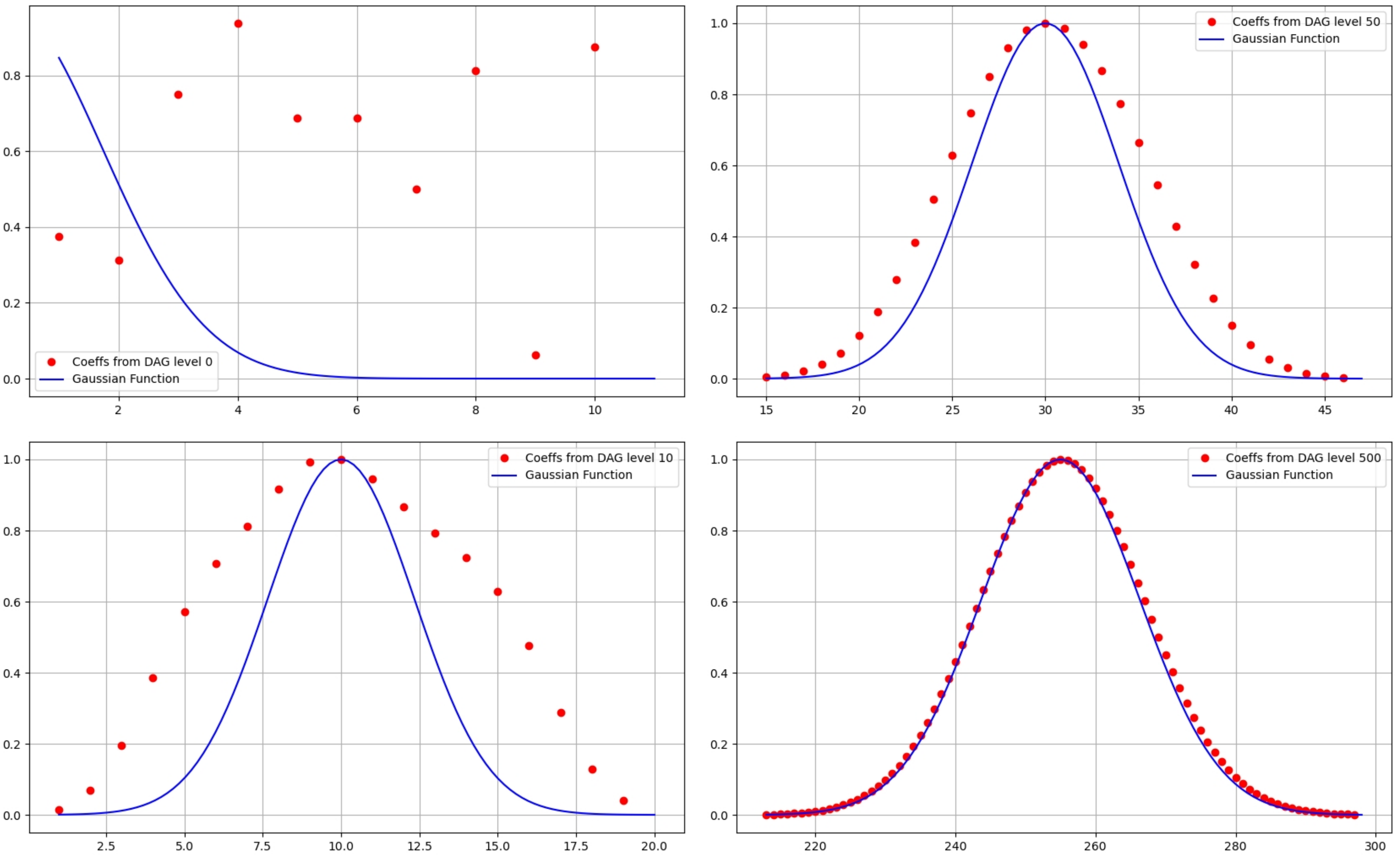}
    \caption{Plot of the normalized node values at level 0, 10, 50, 500. The DAG is initialized with the random basin [16, 6, 5, 12, 15, 11, 11, 8, 13, 1, 14, 7] and node arity of 2}
    \label{fig:dag_vs_gaussian}
\end{figure}

\subsubsection{Additional Invariance Properties}
\begin{figure}[H]
    \centering
    \includegraphics[width=1\linewidth]{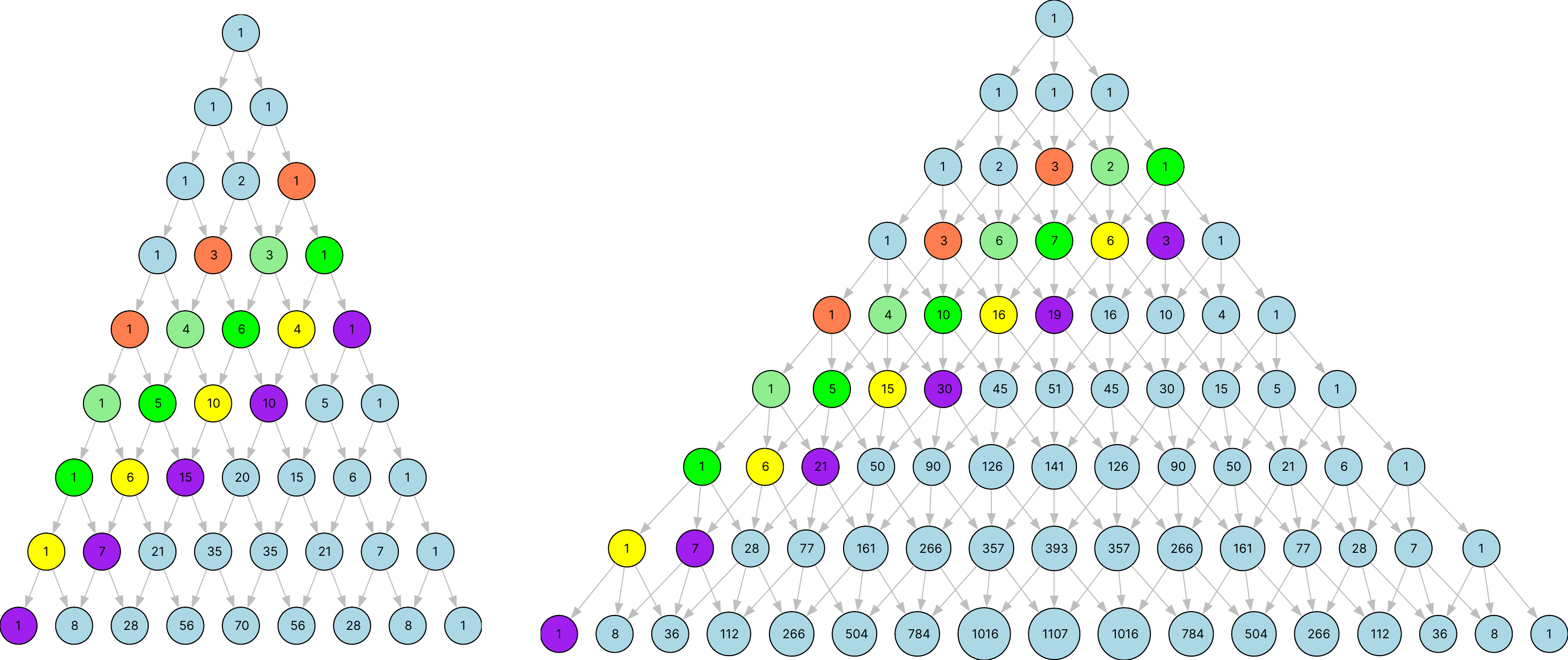}
    \caption{Diagonals of Pascal's triangle on the left. Diagonals of Pascal's ternary triangle on the right.}
    \label{fig:highlighted_diagonals}
\end{figure}

While we do not perform extensive analysis on these apparent DAG Biroot invariance properties, we do note them here. One of the most interesting is that of using diagonal nodes (\cref{fig:highlighted_diagonals}) rather than row nodes to construct the biroot. Our analysis indicates that as larger diagonals are chosen, the convergence error decreases for the biroot. In addition, even when choosing diagonals in $k$ steps rather than what is shown above, we still observe convergence. This is yet more evidence that linearly constructed DAGs form the underlying foundations of rational nth root approximations as the whole structure appears to encode the necessary coefficient distributions.

Table \ref{tab:mean_diagonal_dag_cube_convergence_error} presents statistical analysis for the cube root when using the Diagonal-DAG Biroot over the interval $x \in [0, 10^4]$. With each DAG's level of $m=40$, we generate random basin, node arity, and linear function configurations. For each configuration, we sampled $10^4$ uniformly distributed points and computed the mean absolute error $\mu$ and standard deviation $\sigma$. Here again, we see impressive evidence of convergence.

\begin{remark}
As is shown by Anatriello et al. \cite{Anatriello2022}, the sums of the diagonals of triangles with node arity $k$ form $k$-Padovan-like sequences. For instance, the sum of the diagonals of Pascal's Triangle yields the Fibonacci Sequence. The sum of the diagonals of the ternary triangle yields the Tribonacci sequence following the recurrence relation $T_k = T_{k-1} + T_{k-2} + T_{k-3}$. The pattern continues for higher node arity values $k$.
\end{remark}

\section{Computational Framework}
In order to facilitate both thorough and reproducible analysis of the biroot method, we developed a two tier Python framework. The first framework (MathFlow) is general and aims to provide an ergonomic and comprehensive interface to both symbolic and numerical algorithms in an OOP (Object Orientated Programming) paradigm. We develop it using SymPy (Symbolic CAS), SciPy (for Numerical Algorithms), and NumPy (for underlying Numerical support), with the main goal of unifying them into a single cohesive interface. We also extend these tools with our own custom algorithms such as Padé approximant implementations (Symbolic, Verbose-Symbolic, Numerical), or more advanced root finding methods (hybrid symbolic-numerical). The architectural design goals for this general framework are extensive and deserve their own comprehensive treatment elsewhere. The second framework, however, is an extension of the first. It is what implements the various specific tools for the results described in this paper. The two core contributions we present are (1) code for constructing biroot functions, and (2) code for creating PT-like, but heavily generalized, DAG structures.

\begin{figure}[H]
    \centering
    \includegraphics[width=0.75\linewidth]{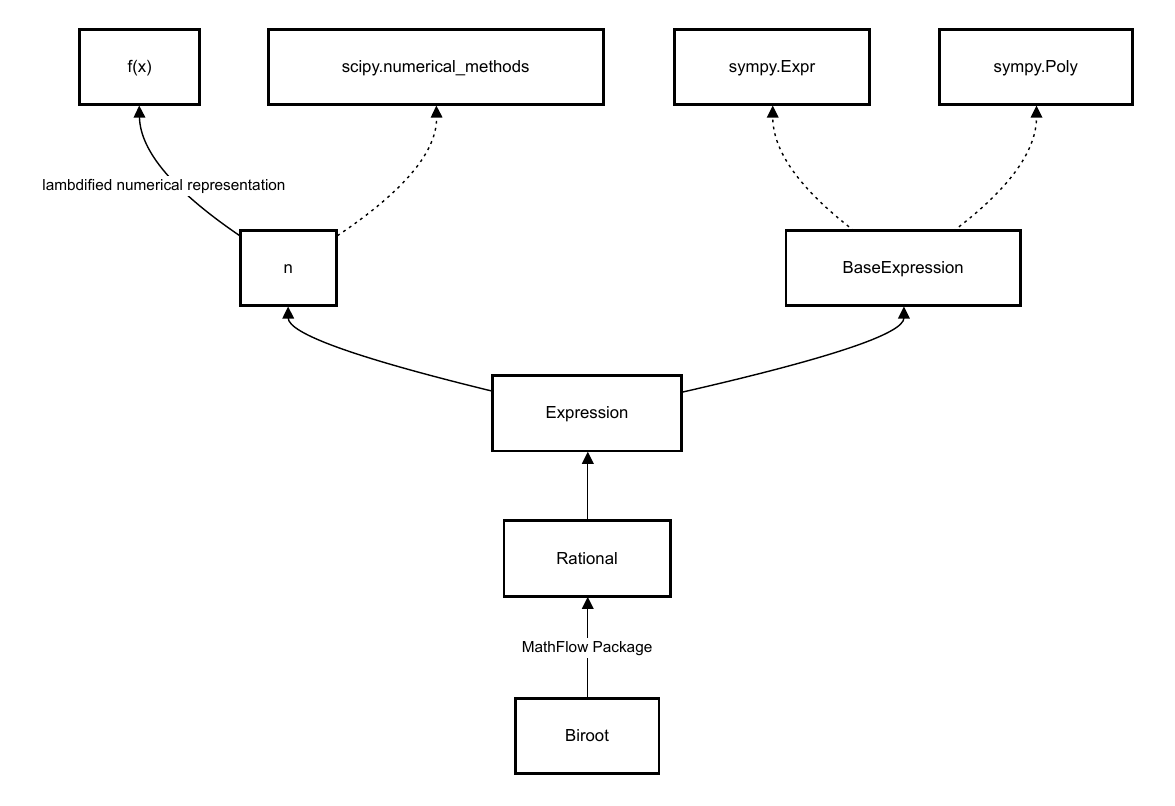}
    \caption{Diagram of the inheritance/proxy nature of the Biroot class. Note that dotted lines mean "Proxy to".}
    \label{fig:biroot_framework_architecture}
\end{figure}

Biroot functions, for instance, can be constructed as objects using
\begin{lstlisting}
>>> from sympy.abc import c  # to parameterize c
>>> b = Biroot(m=12, n=3, c=c)
>>> b.print('latex')
\end{lstlisting}
which yields
\[\frac{c^{12} + 220 c^{9} x + 924 c^{6} x^{2} + 220 c^{3} x^{3} + x^{4}}{12 c^{11} + 495 c^{8} x + 792 c^{5} x^{2} + 66 c^{2} x^{3}}\]
If one then wants to find the error statistics of a certain interval, they can run
\begin{lstlisting}
>>> import numpy as np
>>> e = lambda x: np.abs(b(x) - x**(1/3))  # error function
>>> x = np.linspace(0, 50, 100)  # 100 points over [0, 50]
>>> y = e(x)  # error evaluated over all points x
>>> mean = np.mean(y)
>>> std = np.std(y)
\end{lstlisting}
This works because \verb|b| supports both scalar and vectorized function evaluation. As can be seen in \cref{fig:biroot_framework_architecture}, Biroot objects support all applicable MathFlow (and by extension SymPy and SciPy algorithms). This enables extensive and friction free analysis of the method.

As for generating DAG structures, we created two different data structures: \verb|Level| and \verb|DAG|. The \verb|DAG| object is made up of \verb|Level|s and uses each \verb|Level|'s methods to construct new ones. Each \verb|Level| extends Python's built in \verb|list| class and adds support for continuous function interpolation (to support Gaussian interpolation) and multi-head traversal. A \verb|DAG| instance includes methods for setting continuous interpolation properties, graph construction, level extraction, and diagonal extraction. We have also developed an external module that can be used alongside the \verb|DAG| objects for visualization using the DOT language in Graphviz \footnote{https://graphviz.org/}. For example, to create the DAG
\begin{figure}[H]
    \centering
    \includegraphics[width=0.5\linewidth]{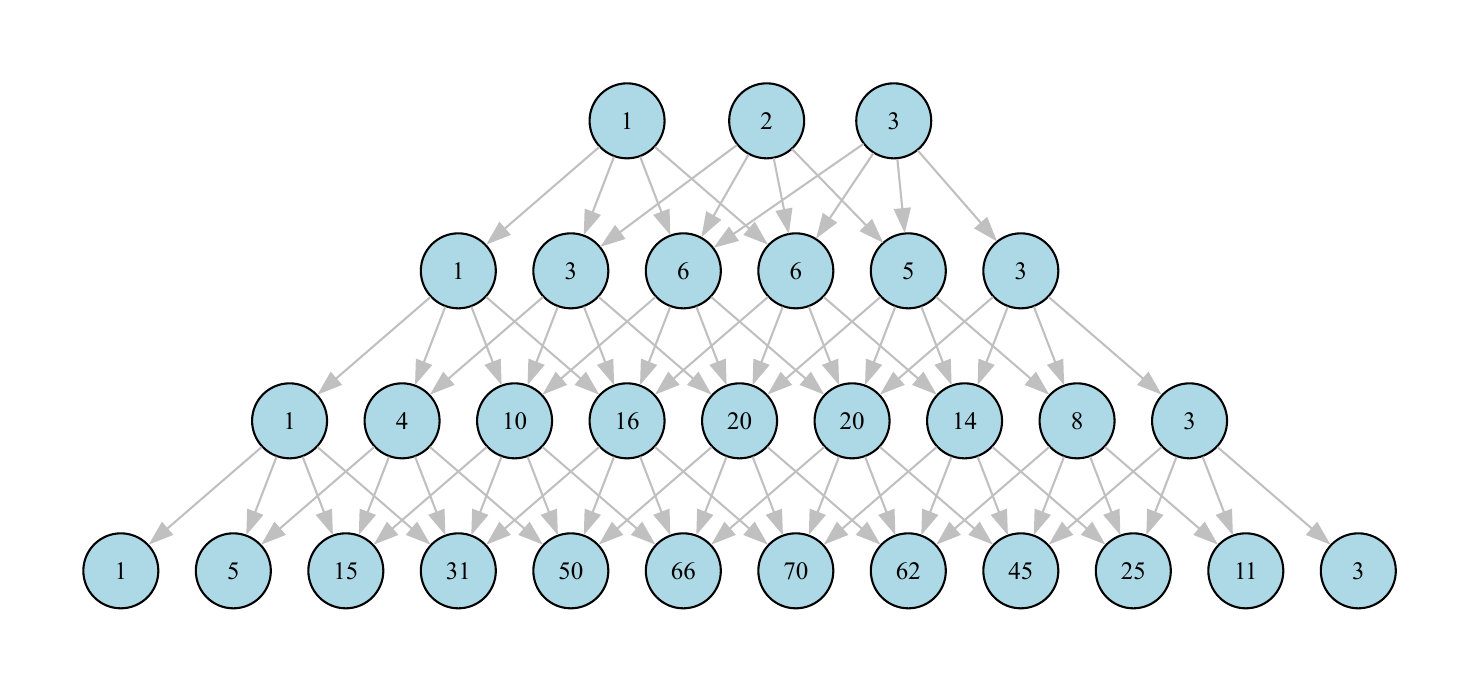}
\end{figure}
\noindent, one would run the code
\begin{lstlisting}
>>> dag = DAG().as_graph(basin=[1, 2, 3], depth=3,
...                      node_arity=4, func=sum)
>>> dot = create_graphviz_from_DAG(dag, node_arity=4)
>>> dot.render(filename='figure', format='pdf')  # save as .pdf
\end{lstlisting}

To create the heatmaps and plots used throughout this paper, we used Seaborn \cite{Waskom2021} and Matplotlib \cite{MDT2025}. After generating multiple datasets containing over 30 million Biroot evaluations across varying intervals and saving them as \verb|.npz| files (NumPy matrices),  we wrote a relatively simple script to create heatmaps. We also saved the best Seaborn configurations so that the figures can automatically be generated and saved upon changes to the data.

All source code, additional materials, and data can be found in the supplementary materials at \href{https://doi.org/10.5281/zenodo.16878055}{https://doi.org/10.5281/zenodo.16878055}.

\section{Applications and Future Work}
Potential applications for the biroot extend beyond just highly accurate nth root formulas. Even less accurate (lower degree) biroot functions appear to still achieve relatively low convergence error over remarkably wide intervals. As a result, they could be used to generate high-quality seed values that require fewer Newton iterations to achieve target precision, potentially reducing overall computational cost. This is a particularly interesting direction to explore in regard to the Gaussian Biroot variant due to its impressive convergence properties. Even a cube biroot function of just degree 14, has a mean absolute convergence error of 0.4105 and standard deviation 0.2242 over a massive interval of zero to one hundred million (it probably extends beyond that, but we did not test to see). Perhaps this could be used to create a hardware accelerator (using Horner's method to optimize the polynomials) for nth root calculations.

Another interesting direction to explore is the Recursive Biroot variant. For instance, $\beta{(x, \beta{(x, \beta{(x, ...)})}}$ could result in accelerated recursive methods similar to newtons method for the nth roots. While we have not tested this much, after some preliminary analysis we do see it making "leaps" towards the correct answer (especially true for large $x$).

Because the biroot method provides closed-form nth root formulas, they could potentially be implemented as custom arithmetic units using ASIC/FPGA designs for dedicated nth root hardware. This approach offers deterministic latency and high parallelization potential, making it particularly suitable for applications requiring massive concurrent nth root computations, such as scientific computing, graphics processing, or machine learning accelerators.

Several theoretical questions merit further investigation. The general nth root Binomial Biroot Conjecture requires rigorous proof, likely employing roots of unity filters and asymptotic analysis techniques. Establishing rigorous and predictable convergence properties across different parameter regimes remains an open problem. The systematic connections between biroot methods and established approximation techniques (particularly Padé approximations, Chebyshev-Padé variants, and Newton's method) need formal characterization and comparison. Most importantly, providing theoretical foundations for the Gaussian Biroot variant would validate its role as an underlying attracting structure for nth root approximation. Finally, the invariance properties observed in DAG structures and their relationship to Gaussian distributions warrant comprehensive theoretical analysis.

\section*{Conclusion}
This paper introduces the biroot method, an approach to nth root approximation using closed-form rational functions with coefficients derived from combinatorial structures. We establish three variants: the Binomial-based Biroot, the Gaussian Biroot leveraging continuous distributions for superior convergence, and the DAG Biroot extending to generalized graph structures. Our theoretical contributions include proof of square root convergence, the Optimal Condition Theorem, and the generalized biroot conjectures supported by extensive computational evidence.

Key open questions include proving the general Binomial Biroot Conjecture and establishing theoretical foundations for the Gaussian variant's exceptional convergence properties. The biroot method offers a promising framework connecting combinatorial mathematics to practical nth root computation.

\section*{Code, Data, and Supplementary Materials}
All source code, additional materials (such as interactive demonstrations), and data can be found in the supplementary materials at \href{https://doi.org/10.5281/zenodo.16878055}{https://doi.org/10.5281/zenodo.16878055}.

\section*{Acknowledgments}
The author wishes to acknowledge the helpful guidance and support from Dr. Kenneth Caviness and other advisors/colleagues at Southern Adventist University. The author also thanks the open-source mathematical software community for the computational tools that enabled this research.

This work is available as a preprint on arXiv (\href{https://arxiv.org/abs/2508.14095}{https://arxiv.org/abs/2508.14095}), and Zenodo (\href{https://doi.org/10.5281/zenodo.16891613}{https://doi.org/10.5281/zenodo.16891613}).

\section*{Declaration of Interest}
The author declares no conflicts of interest.

\bibliographystyle{plainurl}
\bibliography{lib}

\section*{Appendices}
An interesting, but not very useful conjecture (because we don't bother to prove it), is listed due its potential to help in analyzing the error at 1 as the expansion parameter increases for the square root case.

\begin{conjecture}
The factored square biroot of $1$
\[
\beta_m^2{(1,c)} = \frac{(c+1)^m+(c-1)^m}{(c+1)^m-(c-1)^m}
\]
\end{conjecture}

\subsection*{Figures}

\subsection*{Tables}
\begin{table}[H]
    \centering
    \small
    \begin{tabular}{ccc|ccc}
         Basin & Node Arity & Coefficient Vector & Mean Error ($\mu$) & Std. Dev. ($\sigma$)  &Max Error\\
         \hline
         $[7, 3, 9, 7, 5, 2, 2]$ & 3 & $[2, 7, 8]$ & 3.5e-04 & 0.0006 &0.0282\\
         $[3, 16, 1]$ & 4 & $[2, 5, 4, 6]$ & 8.6e-03 & 0.0118 &0.0411\\
         $[6, 7, 5, 4, 2]$ & 4 & $[2, 5, 1, 8]$ & 1.1e-02 & 0.0085 &0.1622\\
         $[14]$ & 2 & $[7, 3]$ & 2.2e+00 & 1.6229 &5.3556\\
         $[2]$ & 3 & $[3, 7, 3]$ & 3.3e-02 & 0.0328 &0.0946\\
         $[6, 8, 4, 8, 4, 4]$ & 4 & $[4, 9, 6, 2]$ & 3.2e-02 & 0.0301 &0.0837\\
         $[1, 1, 14, 5]$ & 4 & $[9, 5, 5, 3]$ & 6.7e+00 & 5.3154 &16.959\\
         $[11, 3, 7, 12, 8]$ & 2 & $[1, 8]$ & 1.9e-05 & 0.0019 &0.1897\\
         $[13]$ & 4 & $[4, 9, 4, 1]$ & 4.3e-02 & 0.0521 &0.1718\\
         $[14]$ & 5 & $[1, 9, 6, 4, 6]$ & 3.8e-06 & 0.0004 &0.0375\\
         $[7, 16]$ & 4 & $[4, 8, 5, 7]$ & 5.1e-02 & 0.0681 &0.2403\\
         $[4, 9, 5, 8, 2]$ & 4 & $[1, 5, 1, 8]$ & 4.8e-05 & 0.0014 &0.1379\\
    \end{tabular}
    \caption{Statistical convergence analysis for DAG cube Biroot approximation over $x \in [0, 10^4]$ with $10^4$ sample points and randomized DAG configurations for $m=40$. Similar results are observed for the fourth root.}
    \label{tab:mean_dag_cube_convergence_error}
\end{table}

\begin{table}[H]
    \centering
    \small
    \begin{tabular}{ccc|ccc}
         Basin & Node Arity & Coefficient Vector & Mean Error ($\mu$) & Std. Dev. ($\sigma$)  & Max Error\\
         \hline
         $[3, 9, 7, 3]$ & 4 & $[8, 5, 7, 2]$ & 8.1e-01 & 0.5681 & 1.6073 \\
         $[9, 3]$ & 2 & $[4, 3]$ & 8.2e-01 & 0.7102 & 2.3360 \\
         $[4, 4, 5, 2]$ & 4 & $[2, 3, 1, 6]$ & 6.0e-03 & 0.0094 & 0.8000 \\
         $[3, 2, 4, 6]$ & 3 & $[5, 4, 8]$ & 2.0e-03 & 0.0047 & 0.3967 \\
         $[2, 1, 8, 7, 1, 3, 8]$ & 2 & $[5, 8]$ & 2.1e-02 & 0.0213 & 0.3216 \\
         $[6]$ & 4 & $[7, 5, 2, 8]$ & 1.3e-02 & 0.0129 & 0.8205 \\
         $[3, 6, 7, 1, 6]$ & 4 & $[4, 5, 8, 7]$ & 2.4e-03 & 0.0030 & 0.1503 \\
         $[7, 4, 4, 3, 5, 9, 1, 9]$ & 4 & $[9, 2, 5, 4]$ & 6.4e-02 & 0.0508 & 0.1329 \\
         $[9, 3, 7]$ & 3 & $[5, 1, 9]$ & 2.8e-02 & 0.0381 & 1.9286 \\
         $[5, 5, 4, 3, 1, 5, 9, 4]$ & 3 & $[3, 7, 9]$ & 2.5e-03 & 0.0036 & 0.2872 \\
         $[9, 4]$ & 3 & $[3, 1, 5]$ & 2.2e-02 & 0.0174 & 0.6065 \\
         $[3, 1, 3, 5]$ & 3 & $[5, 5, 7]$ & 8.2e-03 & 0.0104 & 0.2483 \\
    \end{tabular}
    \caption{Statistical convergence analysis for Diagonal DAG cube biroot over $x \in [0, 10^4]$ with $10^4$ sample points and randomized DAG configurations for $m=40$}
    \label{tab:mean_diagonal_dag_cube_convergence_error}
\end{table}

\end{document}